\documentclass{article}

\usepackage[totalwidth=480pt, totalheight=680pt]{geometry}

\usepackage[leqno]{amsmath}
\usepackage{amsfonts,amssymb,amsthm}
\usepackage{amscd,amsxtra}

\usepackage[all]{xy}
\SelectTips{cm}{} 

\usepackage{eucal}
\usepackage{mathrsfs} 

\newtheorem{theorem}{Theorem}[subsection]
\newtheorem{proposition}[theorem]{Proposition}
\newtheorem{lemma}[theorem]{Lemma}
\newtheorem{corollary}[theorem]{Corollary}

\theoremstyle{definition}
\newtheorem{definition}[theorem]{Definition}
\newtheorem{example}[theorem]{Example}
\theoremstyle{remark}
\newtheorem{remark}[theorem]{Remark}

\numberwithin{equation}{section}

\newcommand{\curly}[1]{\mathscr{#1}}

\newcommand{\field}[1]{\ensuremath{\mathbb{#1}}}
\newcommand{\ZZ}{\field{Z}}
\newcommand{\QQ}{\field{Q}}
\newcommand{\RR}{\field{R}}
\newcommand{\CC}{\field{C}}

\newcommand{\TT}{\field{T}}

\DeclareMathOperator{\Id}{Id}

\DeclareMathOperator{\Tot}{Tot}

\newcommand{\sheafhom}{\underline{\operatorname{Hom}}}

\DeclareMathOperator{\cone}{Cone}
\DeclareMathOperator{\gr}{Gr}
\renewcommand{\Im}{\operatorname{Im}}
\renewcommand{\Re}{\operatorname{Re}}
\DeclareMathOperator{\GL}{GL}

\DeclareMathOperator{\pic}{Pic}

\newcommand{\sheafaut}{\underline{\operatorname{Aut}}}

\newcommand{\gerbe}[1]{\mathcal{#1}}
\newcommand{\twogerbe}[1]{\mathsf{#1}}

\newcommand{\herm}[1]{\underline{\operatorname{herm}}({#1})}
\newcommand{\conn}[1]{\underline{\operatorname{Co}}({#1})}

\newcommand{\aut}[1]{\sheafaut (#1)}

\newcommand{\tors}[1]{\operatorname{\mathsf{Tors}}(#1)}
\newcommand{\gerbes}[1]{\operatorname{\mathsf{Gerbes}}(#1)}

\newcommand{\sheaf}[1]{\underline{\mathnormal{#1}}}
\newcommand{\sha}[2][\bullet]{\sheaf{A}_{#2}^{#1}}
\newcommand{\she}[2][\bullet]{\sheaf{\mathcal{E}}_{#2}^{#1}}
\newcommand{\sho}[1]{\mathcal{O}_{#1}}
\newcommand{\shomega}[2][\bullet]{\sheaf{\Omega}_{#2}^{#1}}
\newcommand{\deligne}[3][\bullet]{#2(#3)^{#1}_\mathcal{D}}
\newcommand{\deltilde}[3][\bullet]{%
  \smash[t]{\widetilde{#2(#3)}}^{#1}_\mathcal{D}}

\newcommand{\dhh}[2][\bullet]{D(#2)_\mathit{h.h.}^{#1}}
\newcommand{\delH}[4][\bullet]{H^{#1}_\mathcal{D}(#2, #3(#4))}
\newcommand{\dhhH}[3][\bullet]{%
  H^{#1}_{\mathcal{D}_{\smash[b]{\mathit{h.h.}}}}(#2,#3)}

\newcommand{\qi}{\overset{\simeq}{\rightarrow}}
\newcommand{\lqi}{\overset{\simeq}{\longrightarrow}}

\newcommand{\hyper}[1]{\mathbf{#1}}
\newcommand{\HHH}{\hyper{H}}

\newcommand{\cover}[1]{\mathfrak{#1}}
\newcommand{\vC}[3][\bullet]{\Check{C}^{#1}(#2,#3)}

\newcommand{\complex}[1]{\mathsf{#1}} 
\newcommand{\CCC}{\complex{C}}

\newcommand{\del}{\partial} 
\newcommand{\delb}{\Bar\partial} 
\newcommand{\deltacheck}{\delta} 

\newcommand{\hodge}[1]{\curly{#1}}
\newcommand{\hM}{\hodge{M}}
\newcommand{\hH}{\hodge{H}}

\newcommand{\lto}{\longrightarrow}
\newcommand{\iso}{\cong} 
\newcommand{\coin}{\equiv} 
\newcommand{\unit}{\times} 
\newcommand{\tame}[2]{\bigl(#1,#2\bigr]}
\newcommand{\tamehh}[2]{\bigl(#1,#2\bigr]_{\smash[b]{\mathit{h.h.}}}}
\newcommand{\tameg}[2]{\bigl\lbrace #1,#2\bigr\rbrace}
\newcommand{\tate}{2\pi\smash{\sqrt{-1}}}

\newcommand{\abs}[1]{\left\lvert#1\right\rvert}
\newcommand{\norm}[1]{\lVert#1\rVert}

\newcommand{\onehalf}{\frac{1}{2}}

\newcommand{\eqdef}{\:\overset{\mathrm{def}}{=}\:}
\newcommand{\li}{\mathrm{Li}_2} 

\newcommand{\bei}{Be\u\i{}linson}
\newcommand{\cech}{\v{C}ech}

\usepackage{ifpdf}
\usepackage{ifthen}
\ifthenelse{\boolean{pdf}}{%
  \usepackage[%
              hyperindex=false,
              plainpages=false,
              colorlinks=true,
              hypertexnames=true,
              plainpages=false,
              backref=page,
              pdfpagelabels]{hyperref}}{%
  \usepackage[hypertex]{hyperref}}

\title{Hermitian-holomorphic ($2$)-Gerbes and tame symbols}

\author{Ettore Aldrovandi\\
  Department of Mathematics\\
  Florida State University\\
  Tallahassee, FL 32306-4510, USA\\
  \texttt{aldrovandi@math.fsu.edu} }

\date{}

\begin{document}

\maketitle

\begin{abstract}
  The tame symbol of two invertible holomorphic functions can be
  obtained by computing their cup product in Deligne cohomology,
  and it is geometrically interpreted as a holomorphic line
  bundle with connection. In a similar vein, certain higher tame
  symbols later considered by Brylinski and McLaughlin are
  geometrically interpreted as holomorphic gerbes and $2$-gerbes
  with abelian band and a suitable connective structure.
  
  In this paper we observe that the line bundle associated to the
  tame symbol of two invertible holomorphic functions also
  carries a fairly canonical hermitian metric, hence it
  represents a class in a Hermitian holomorphic Deligne
  cohomology group.
  
  We put forward an alternative definition of hermitian
  holomorphic structure on a gerbe which is closer to the
  familiar one for line bundles and does not rely on an explicit
  ``reduction of the structure group.'' Analogously to the case
  of holomorphic line bundles, a uniqueness property for the
  connective structure compatible with the hermitian-holomorphic
  structure on a gerbe is also proven. Similar results are
  proved for $2$-gerbes as well.
  
  We then show the hermitian structures so defined propagate to a
  class of higher tame symbols previously considered by Brylinski
  and McLaughlin, which are thus found to carry corresponding
  hermitian-holomorphic structures. Therefore we obtain an
  alternative characterization for certain higher Hermitian
  holomorphic Deligne cohomology groups.
\end{abstract}

\tableofcontents

\section{Introduction}
\label{sec:Introduction}

The aim of this work is two-fold. For an analytic manifold $X$ we
investigate geometric objects corresponding to the elements of
certain low-degree Hermitian-Holomorphic Deligne cohomology
groups. These groups, denoted here $\dhhH[k]{X}{l}$, for two
integers $k$ and $l$, were defined in~\cite{bry:quillen} and, in
a slightly different fashion, later in~\cite{math.CV/0211055}. It
is already an observation by Deligne (cf. \cite{esn:char}) that
\begin{math}
  \dhhH[2]{X}{1} \iso \widehat{\pic X}\,,
\end{math}
the group of isomorphism classes of holomorphic line bundles with
hermitian fiber metric. Here we define an appropriate notion of
hermitian structure on a gerbe (or $2$-gerbe) bound by
$\sho{X}^\unit$ and show that the corresponding (equivalence)
classes are in bijective correspondence with the elements of
$\dhhH[k]{X}{1}$, for $k=3,4$.

As a second result and application, we show that the torsors and
($2$-)gerbes underlying the cup products in ordinary Deligne
cohomology studied by
Brylinski-McLaughlin~\cite{brymcl:deg4_I,brymcl:deg4_II} can be
equipped in a rather natural way with the above mentioned
hermitian structures, thus producing classes in the
Hermitian-Holomorphic variant. More precisely, we modify the cup
product at the level of Deligne complexes to land into a
Hermitian-Holomorphic one. This modification is actually quite a
natural one from the point of view of Mixed Hodge Structures.

\subsection{Background notions}
\label{sec:background-notions}

To explain things a little bit more, let $X$ be an analytic
manifold and let $A\subseteq \RR$ be a subring---typically
$A=\ZZ,\QQ$ or $\RR$. For any integer $j$, set $A(j)= (\tate)^j
A$ and let $\deligne{A}{j}$ be the Deligne complex
\begin{displaymath}
  A(j)_X \hookrightarrow \sho{X} \to \shomega[1]{X}\to\dots \to
  \shomega[j-1]{X} \,.
\end{displaymath}
It is well known that (at the level of the derived category)
there are maps
\begin{math}
  \deligne{A}{j}\otimes \deligne{A}{k} \to \deligne{A}{j+k}
\end{math}
inducing a cup product in cohomology
\begin{displaymath}
  \delH[p]{X}{A}{j} \otimes \delH[q]{X}{A}{k}
  \xrightarrow{\cup} \delH[p+q]{X}{A}{j+k}\,,
\end{displaymath}
where we have used the notation
\begin{math}
  \delH[p]{X}{A}{j} = \HHH^p (X, \deligne{A}{j} )
\end{math}
for the \emph{Deligne cohomology} groups, and $\HHH^\bullet(X,-)$
denotes hypercohomology.

The question of obtaining a geometric picture of the cup product
in cohomology is a very interesting one. A chief foundational
example is the following. For $A=\ZZ$ the product
\begin{equation}
  \label{eq:54}
  \deligne{\ZZ}{1}\otimes \deligne{\ZZ}{1}
  \lto \deligne{\ZZ}{2}
\end{equation}
corresponds to the morphism
\begin{equation}
  \label{eq:53}
  \sho{X}^\unit \otimes \sho{X}^\unit \lto
  \bigl( \sho{X}^\unit \xrightarrow{d\log} \shomega[1]{X}
  \bigr)
\end{equation}
via the quasi-isomorphisms
\begin{math}
  \deligne{\ZZ}{1} \qi \sho{X}^\unit [-1]
\end{math}
and
\begin{math}
  \deligne{\ZZ}{2} \qi
  \bigl( \sho{X}^\unit \xrightarrow{d\log}
  \shomega[1]{X} \bigr)[-1]\,.
\end{math}
Deligne gave a geometric construction of~\eqref{eq:53}and the
ensuing cup product
\begin{displaymath}
  \sho{X}^\unit (X) \otimes \sho{X}^\unit (X)
  \xrightarrow{\cup}
  \HHH^1 \bigl( X,\sho{X}^\unit\xrightarrow{d\log} \shomega[1]{X}
  \bigr) 
\end{displaymath}
in his work on tame symbols, cf.~\cite{del:symbole}: If $f$ and
$g$ are two invertible functions on $X$, namely two elements of
$\sho{X}^\unit$, their cup product corresponds to a
$\sho{X}^\unit$-torsor, denoted $\tame{f}{g}$, equipped with an
analytic connection. Furthermore, if $X$ is a Riemann surface,
the complex
\begin{math}
  \bigl( \sho{X}^\unit \xrightarrow{d\log} \shomega[1]{X}
  \bigr)
\end{math}
is quasi-isomorphic to $\CC^\unit$ and the product is interpreted
as the \emph{holonomy} of the connection. For $X$ equal to a
punctured disk $D_p$ centered at $p$, if $f$ and $g$ are
holomorphic on $D_p$, meromorphic at $p$, the holonomy of
$\tame{f}{g}$ computes the \emph{tame symbol}
\begin{displaymath}
  (f,g)_p = (-)^{v(f)v(g)}\bigl( f^{v(g)}/g^{v(f)}\bigr) (p)\,,
\end{displaymath}
where $v(f)$ is the valuation of $f$ at $p$,
cf.~\cite{MR86h:11103,del:symbole,MR94k:19002}. This justifies
the use of the name \emph{tame symbol} for $\tame{f}{g}$.

A particularly pleasant property is that when $f$ and $1-f$ are
both invertible a calculation \cite{del:symbole} using the
classical Euler's dilogarithm $\li$ shows that $\tame{f}{1-f}$ is
isomorphic to the trivial torsor equipped with the trivial
connection $d$, namely the unit element in the group
\begin{math}
  \HHH^1 \bigl(X,\sho{X}^\unit \xrightarrow{d\log} \shomega[1]{X}
  \bigr)\,.
\end{math}
{}From this one also builds an interpretation of the symbol
associated to $f$ and $g$ in terms of Mixed Hodge
Structures~\cite{del:symbole}.

In this particular example there appear degree $1$ and $2$
Deligne cohomology groups: specifically, it is made use of the
fact that invertible functions determine elements in the group
\begin{math}
  \delH[1]{X}{\ZZ}{1}\iso \sho{X}^\unit(X)\,,
\end{math}
and, given $f$ and $g$, the class of the torsor with connection
$\tame{f}{g}$ is an element of
\begin{math}
  \delH[2]{X}{\ZZ}{2}\iso
  \HHH^1 \bigl(X,\sho{X}^\unit \xrightarrow{d\log}
  \shomega[1]{X} \bigr)\,.
\end{math}
It is therefore natural to investigate the geometric objects
corresponding to similar cup products of higher degree. The case
of $\tame{f}{L}$, where $f$ is again an invertible function and
$L$ is an $\sho{X}^\unit$-torsor, so it determines a class in 
\begin{math}
  \delH[2]{X}{\ZZ}{1} \iso H^1(X,\sho{X}^\unit)\,,
\end{math}
was already considered in ref.~\cite{del:symbole}, where it is
interpreted in terms of a gerbe $\gerbe{G}$ over $X$.

This idea has been further pursued by
Brylinski-McLaughlin,~\cite{brymcl:deg4_I,brymcl:deg4_II}. In
their study of degree $4$ characteristic classes they considered
the symbols $\tame{f}{L}\in \delH[3]{X}{\ZZ}{2}$ and, for a pair
of $\sho{X}^\unit$-torsors, $\tame{L}{L'}\in
\delH[4]{X}{\ZZ}{2}$.  The corresponding geometric objects are
identified with a gerbe (resp. a $2$-gerbe) both equipped with
the appropriate analog of a connection. Furthermore, the obvious
map 
\begin{math}
  \deligne{\ZZ}{2}\to \deligne{\ZZ}{1}
\end{math}
induces a corresponding map
\begin{math}
  \delH[k]{X}{\ZZ}{2}\to \delH[k]{X}{\ZZ}{1}
\end{math}
which simply forgets the connection. Therefore elements in the
groups
\begin{math}
  \delH[k]{X}{\ZZ}{1}\,,
\end{math}
for $k=3,4$ correspond to equivalence classes of ($2$-)gerbes
bound by $\sho{X}^\unit$,
cf.\cite{MR95m:18006,brymcl:deg4_I,brymcl:deg4_II}. Thus in the
end several Deligne cohomology groups have a concrete
interpretation in terms of geometric data.

\emph{Hermitian-Holomorphic} Deligne cohomology, as defined by
Brylinski, cf.~\cite{bry:quillen}, is an enhanced version of
Deligne cohomology. For all positive integers $l$ Brylinski
introduces certain complexes $C(l)^\bullet$, and defines the
Hermitian-Holomorphic Deligne cohomology groups as the sheaf
hypercohomology groups:
\begin{math}
  \dhhH[k]{X}{l} = \HHH^k(X,C(l)^\bullet)\,.
\end{math}
The complex $C(l)^\bullet$ has a map
\begin{math}
  C(l)^\bullet \to \deligne{\ZZ}{l}\,,
\end{math}
thus there is an obvious map
\begin{math}
  \dhhH[k]{X}{l} \to \delH[k]{X}{\ZZ}{l}
\end{math}
forgetting the extra-structure.

A primary example is provided by Deligne's observation mentioned
before, cf.~\cite{esn:char}, that
\begin{equation}
  \label{eq:55}
  \widehat{\pic X}\iso
  \HHH^2 \bigl( X, \ZZ(1)_X\to \sho{X} \to \she[0]{X} \bigr)\,,
\end{equation}
where $\widehat{\pic X}$ is the set of isomorphism classes of
$\sho{X}^\unit$-torsors with hermitian metric, and $\she[0]{X}$
is the sheaf of smooth real-valued functions on $X$. The complex
in~\eqref{eq:55} is quasi-isomorphic to $C(1)^\bullet$, therefore
\begin{displaymath}
  \widehat{\pic X} \iso\dhhH[2]{X}{1}\,.
\end{displaymath}
In fact, both complexes are quasi-isomorphic to the complex
\begin{math}
  \bigl(\sho{X}^\unit\oplus \sheaf{\TT}_X \to
  \sheaf{\CC}^\unit_X\bigr)[-1]\,,
\end{math}
~\cite{brymcl:deg4_II,bry:quillen}, which encodes the reduction
of the torsor structure from $\sho{X}^\unit$ to $\sheaf{\TT}_X$
afforded by the hermitian metric.

Concerning higher degrees, Brylinski-McLaughlin
\cite{brymcl:deg4_II,MR97d:32041} gave a geometric interpretation
for some of the groups $\dhhH[k]{X}{l}$, $k=3,4$ and $l=1,2$ in
terms of classes of gerbes and $2$-gerbes bound by $\sheaf{\TT}_X$
and equipped with a concept of connection valued in an
appropriate Hodge filtration of the de~Rham complex of $X$.

\subsection{Statement of the results}
\label{sec:statement-results}

In this work we take on the same question of a geometric
interpretation for some Hermitian-Holomorphic Deligne cohomology
groups from a holomorphic view-point which, we believe, is
complementary to that of Brylinski-McLaughlin. We define a
hermitian structure on a $\sho{X}^\unit$-gerbe $\gerbe{G}$ as the
assignment of a $\she[0]{U,+}$-torsor (the ``$+$'' denotes
positive functions) to any object $P$ of $\gerbe{G}_U$ subject to
several conditions spelled out in Definition~\ref{def:2}. We
prove that classes of gerbes with hermitian structures in this
sense correspond to elements of
\begin{math}
  \dhhH[3]{X}{1} \iso
  \HHH^3 \bigl( X, \ZZ(1)_X\to \sho{X} \to \she[0]{X} \bigr)\,,
\end{math}
in complete analogy with~\eqref{eq:55}. Moreover we can define a
type $(1,0)$-connective structure on $\gerbe{G}$ by requiring
that to any object $P$ of $\gerbe{G}_U$ be assigned a
$F^1\!\sha[1]{U}$-torsor, essentially repeating the steps in
ref.~\cite{brymcl:deg4_II}. (Here $\sha{U}$ is the smooth
$\CC$-valued de~Rham complex, and $F^1$ is the first Hodge
filtration.) Then a notion of compatibility between the hermitian
structure and the connective one is defined, and in fact we prove
there is only one such type $(1,0)$ connective structure
compatible with a given hermitian structure, up to equivalence.
This result is analogous to the corresponding statement for
hermitian holomorphic line bundles, that there is a unique
connection --- the \emph{canonical or Griffiths connection} ---
compatible with both structures.

Similar results are available for $2$-gerbes: we define a
hermitian structure for a $\sho{X}^\unit$-$2$-gerbe
$\twogerbe{G}$ as the assignment of a $\she[0]{U,+}$-\emph{gerbe}
for each object $P$ of $\twogerbe{G}_U$, subject to several
conditions spelled out in Definition~\ref{def:5}. Analogously to
the simpler case of gerbes, we have a concept of type $(1,0)$
connectivity compatible with the hermitian structure and a
uniqueness result up to equivalence.

A second line of results is more specific to the tame symbols we
encountered before. Alongside with the map of complexes
\begin{displaymath}
  \deligne{\ZZ}{1} \otimes \deligne{\ZZ}{1} \lto \deligne{\ZZ}{2}
\end{displaymath}
we define a companion map
\begin{equation}
  \label{eq:56}
  \deligne{\ZZ}{1} \otimes \deligne{\ZZ}{1} \lto \tate\otimes
  C(1)^\bullet 
\end{equation}
so that it is possible to obtain a different cup product valued
in Hermitian-Holomorphic Deligne cohomology:
\begin{displaymath}
  \delH[i]{X}{\ZZ}{1} \otimes \delH[j]{X}{\ZZ}{1}
  \overset{\cup}{\lto} \tate\otimes \dhhH[i+j]{X}{1}\,.
\end{displaymath}
An immediate consequence is that for $f$ and $g$ invertible, and
$L,L'$ line bundles, the torsor $\tame{f}{g}$ and the gerbe
$\tame{f}{L}$ support natural hermitian structures of the type
discussed above, in addition to the analytic connection (or
connective) ones associated with the cup product in standard
Deligne cohomology. The same conclusions are valid for the
$2$-gerbe $\tame{L}{L'}$. It turns out that supporting both
structures is an easy consequence of the commutativity of the
following diagram:
\begin{displaymath}
  \begin{CD}
    \delH[i]{X}{\ZZ}{1} \otimes \delH[j]{X}{\ZZ}{1}
    @>{\cup}>>
    \tate\otimes \dhhH[i+j]{X}{1}\\
    @V{\cup}VV  @VV{\mathit{forget}}V\\
    \delH[i+j]{X}{\ZZ}{2}
    @>{\mathit{forget}}>>
    \delH[i+j]{X}{\ZZ}{1}
  \end{CD}
\end{displaymath}
Indeed, forgetting either structure, brings us back to the same
underlying object.

The map~\eqref{eq:56} has a rather natural definition from the
point of view of Mixed Hodge Structures, whose role in the matter
was mentioned in relation with the product~\eqref{eq:54},
see~\cite{del:symbole}. Namely, there is a ``universal'' MHS
$\hM^{(2)}$ corresponding to an iterated extension of $\ZZ(0)$ by
$\ZZ(1)$ by $\ZZ(2)$, where in this case $\ZZ(n)$ denotes a
Hodge-Tate structure. To $\hM^{(2)}$ we can associate a tensor
--- the ``big period'' ---
\begin{math}
  P(\hM^{(2)})\in \CC\otimes_\QQ\CC\,,
\end{math}
cf.~\cite{MR99i:19004}. The period is in fact a multiple of the
extension class of $\hM^{(2)}$, and it belongs to the kernel
\begin{math}
  \curly{I} = \ker\bigl(
  m\colon \CC \otimes_\QQ \CC \to \CC
  \bigr)
\end{math}
of the multiplication map. We find the map~\eqref{eq:56}
corresponds to the image of $P(\hM^{(2)})$ under the ``imaginary
part'' projection 
\begin{math}
  \CC\otimes_\QQ \CC \to \RR (1)
\end{math}
given by $a\otimes b\mapsto \Im (a) \Re (b)$. On the other hand,
the standard one~\eqref{eq:54} involves the projection onto the
K\"ahler differentials 
\begin{math}
  \curly{I} \to \curly{I}/\curly{I}^2
\end{math}
given by $a\otimes b \mapsto a\, db$.

Another consequence of the previous diagram is that
$\tame{f}{g}$, $\tame{f}{L}$, and $\tame{L}{L'}$ come equipped
with two connection (or connective) structures. If the unitary
connection in a line bundle $L$ is also analytic, then $L$ is
flat. In the case of $\tame{f}{g}$ we find there is an
obstruction to this type of compatibility. This can be cast in
cohomological terms, which allows to extend these considerations
to $\sho{X}^\unit$-gerbes and $2$-gerbes. We find that the
obstruction vanish, so compatibility can always be achieved.

\subsection{Outline of the paper}
\label{sec:organization}

This work is organized as follows. In section~\ref{sec:Prelim} we
make some preliminaries observations about Deligne complexes and
cohomology and collect a few needed facts. We recall the
definition of Hermitian-Holomorphic Deligne cohomology and state
some of its properties in section~\ref{sec:herm-holom-deligne}.
Alongside Brylinski's complex $C(l)^\bullet$, we use a complex
quasi-isomorphic to it, denoted $\dhh{l}$, which for a line
bundle directly encodes the data defining the \emph{canonical
  connection.}

In section~\ref{sec:tame-symb-herm} we recall the definition of
the tame symbol $\tame{f}{g}$ for two invertible functions and
some of its properties. We define the modified
product~\eqref{eq:56} and show that through it, the torsor
associated to $\tame{f}{g}$ also comes equipped with a hermitian
structure. As mentioned before, the product~\eqref{eq:56} and its
relation with the standard for Deligne complexes become more
clear when analyzed in terms of Hodge Structures. In order to do
this, we felt necessary to recall a few elementary facts and
calculations concerning Hodge-Tate structures that are certainly
well-known to experts. For this reason, and also because this
development lies somewhat aside this work's main lines, we
present this material in appendix~\ref{sec:remarks-hodge-tate}.
This presentation relies in part on the Heisenberg group picture
of the Deligne torsor, which we have recalled in
section~\ref{sec:heisenberg-group}.

Section~\ref{sec:herm-holorm-gerb} is the main part of this work.
There we redefine the notion of hermitian structure (modeled
after that of connective structure) and prove that equivalence
classes of these are classified by the groups $\dhhH[k]{X}{1}$.
We then apply this classification to the Hermitian structures and
the product~\eqref{eq:56} for the higher versions of the tame
symbols considered by Brylinski-McLaughlin.

The interplay between the analytic connection (or connective)
structures arising from standard Deligne cohomology and their
hermitian counterparts defined here is analyzed in
sections~\ref{sec:comparisons} and~\ref{sec:comp-relat-with}.

\subsection*{Acknowledgments}

Parts of the present work were written while visiting the
Department of Mathematics, Aarhus University, \AA rhus, Denmark;
the International School for Advanced Studies (SISSA), Trieste,
Italy; the Department of Mathematics, Instituto Superior
T\'ecnico, Lisbon, Portugal. It is a pleasure to thank all these
institutions for hospitality, support, and for providing an
excellent, friendly, and stimulating research environment. It is
also a pleasure to thank the anonymous referee for raising
important points and providing several stimulating comments
leading to a much improved version of the paper.

\section{Preliminaries}
\label{sec:Prelim}

\subsection{Notation and conventions}
\label{sec:Notation-conventions}

If $z$ is a complex number, then $\pi\sb{p} (z) \eqdef \onehalf (
z + (-1)\sp p \Bar z)$, and similarly for any other complex
quantity, e.g.  complex valued differential forms. For a subring
$A$ of $\RR$ and an integer $j$, $A(j) = (\tate)^j\,A$ is the
Tate twist of $A$. We identify $\CC/\ZZ(j) \iso \CC^\unit$ via
the exponential map $z \mapsto \exp (z/(\tate)^{j-1})$, and $\CC
/ \RR(j) \iso \RR (j-1)$.

If $X$ is a complex manifold, $\sha{X}$ and $\shomega{X}$ denote
the de~Rham complexes of sheaves of smooth $\CC$-valued and
holomorphic forms, respectively. We denote by $\she{X}$ the
de~Rham complex of sheaves of real valued differential forms and
by $\she{X} (j)$ the twist $\she{X} \otimes_\RR \RR(j)$. We set
$\sho{X} \coin \shomega[0]{X}$ as usual.  When needed,
$\sha[{p,q}]{X}$ will denote the sheaf of smooth $(p,q)$-forms.
We use the standard decomposition $d=\del + \delb$ according to
types. Furthermore, we introduce the differential operator $d^c =
\del -\delb$ (contrary to the convention, we omit the factor
$1/(4\pi \sqrt{-1})$). We have $2\del\delb = d^cd$. The operator
$d^c$ is an imaginary one and accordingly we have the rules
\begin{displaymath}
  d\pi_p(\omega) = \pi_p(d\omega)\,,\quad
  d^c\pi_p(\omega) = \pi_{p+1}(d^c\omega)
\end{displaymath}
for any complex form $\omega$.

An open cover of $X$ will be denoted by $\cover{U}_X$. If
$\{U_i\}_{i\in I}$ is the corresponding collection of open sets,
we write $U_{ij} = U_i\cap U_j$, $U_{ijk} = U_i\cap U_j\cap
U_k$, and so on. More generally we can also have $\cover{U}_X =
\{ U_i \to X\}_{i\in I}$, where the maps are regular coverings in
an appropriate category. In this case intersections are replaced
by $(n+1)$-fold fibered products
\begin{math}
  U_{i_0 i_1\dotsb i_n}
  = U_{i_0} \times_X\dotsb \times_X U_{i_n}\,.
\end{math}

If $\sheaf{F}^\bullet$ is a complex of abelian sheaves on $X$,
its \cech\ resolution with respect to a covering $\cover{U}_X\to
X$ is the double complex
\begin{displaymath}
  \CCC^{p,q} (\sheaf{F}) \eqdef
  \vC[q]{\cover{U}_X}{\sheaf{F}^p}\,,
\end{displaymath}
where the $q$-cochains with values in $\sheaf{F}^p$ are given by
\begin{math}
  \prod \sheaf{F}^p (U_{i_0\dotsb i_n})\,.
\end{math}
The \cech\ coboundary operator is denoted $\deltacheck$. The
convention we use is to put the index along the \cech\ resolution
in the \emph{second} place, so if we denote by $d$ the
differential in the complex $\sheaf{F}^\bullet$, the total
differential is given by $D=d + (-1)^p \deltacheck$ on the
component $\vC[q]{\cover{U}_X}{\sheaf{F}^p}$ of the total
simple complex. Furthermore, recall that the Koszul sign rule
causes a sign being picked whenever two degree indices are
formally exchanged. For \cech\ resolutions of complexes of
sheaves it leads to the following conventions. If
$\sheaf{G}^\bullet$ is a second complex of sheaves on $X$, then
one defines the cup product
\begin{displaymath}
  \cup : \CCC^{p,q}(\sheaf{F}) \otimes \CCC^{r,s}(\sheaf{G})
  \lto
  \vC[q+s]{\cover{U}_X}{\sheaf{F}^p\otimes \sheaf{G}^r} \subset
  \CCC^{p+r,q+s}(\sheaf{F}\otimes\sheaf{G})
\end{displaymath}
of two elements $\{f_{i_0,\dotsc,i_q}\}\in \CCC^{p,q}(\sheaf{F})$
and $\{g_{j_0,\dotsc,j_s}\} \in \CCC^{r,s}(\sheaf{G})$ by
\begin{displaymath}
  (-1)^{qr}\,f_{i_0,\dots,i_q}\otimes
  g_{i_q,i_{q+1},\dotsc,i_{q+s}} \,.
\end{displaymath}
For a given complex of abelian objects, say $\CCC^\bullet$, the
symbol $\sigma^{i}$ denotes sharp truncation at the index $i$:
$\sigma^i\CCC^p=0$ for $p<i$.

\subsection{Deligne cohomology}
\label{sec:deligne-cohomology}

There are several models for the complexes to use to compute
Deligne cohomology \cite{esn-vie:del,MR86h:11103}. For $A\subset
\RR$ and an integer $j$ the latter is the hypercohomology:
\begin{displaymath}
  \delH{X}{A}{j} = \HHH^\bullet (X, \deligne{A}{j})\,. 
\end{displaymath}
Here $\deligne{A}{p}$ is the Deligne complex
\begin{align}
  \label{eq:1}
  \deligne{A}{j} &=  A(j)_X
  \overset{\imath}{\lto} \sho{X}
  \overset{d}{\lto} \shomega[1]{X}
  \overset{d}{\lto} \dotsm
  \overset{d}{\lto} \shomega[{j-1}]{X} \\
  \label{eq:2}
  &\lqi  \cone \big( A(j)_X \oplus F^j\shomega{X}
  \xrightarrow{\imath -\jmath} \shomega{X} \big) [-1]\,,
\end{align}
where $F^j\shomega{X}$ in eqn.~\eqref{eq:2} is the Hodge
(``stupid'') filtration on the de~Rham complex. The symbol $\lqi$
denotes a quasi-isomorphism. In view of \bei\ formula for the cup
product on cones to be recalled below~\cite{bei:hodge_coho},
Deligne complexes acquire a family of cup-products (depending on
a real parameter $\alpha$)
\begin{displaymath}
  \deligne{A}{j} \otimes \deligne{A}{k}
  \overset{\cup_\alpha}{\lto}
  \deligne{A}{j+k}\,.
\end{displaymath}
Cup products related to different values of the parameter
$\alpha$ are related by homotopy-commutative diagrams, hence they
induce a well defined graded commutative cup-product in
cohomology
\begin{equation}
  \label{eq:3}
    \delH[p]{X}{A}{j} \otimes \delH[q]{X}{A}{k}
  \xrightarrow{\cup} \delH[p+q]{X}{A}{j+k}\,.
\end{equation}
In order to explicitly compute cup products, the model given by
eq.~\eqref{eq:1} leads to simpler formulas (when it can be
used). If $f\in \deligne{A}{j}$ and $g\in \deligne{A}{k}$, then
from ref.~\cite{esn-vie:del} we quote:
\begin{equation}
  \label{eq:4}
  f\cup g =
  \begin{cases}
    f\cdot g & \deg f = 0\,,\\
    f\wedge dg & \deg f > 0\;
    \text{and}\; \deg g =k\,,\\
    0 &\text{otherwise.}
  \end{cases}
\end{equation}
The following examples are well known and will frequently recur
in the following.
\begin{example}
  For $A=\ZZ$ it is immediately verified that
  \begin{math}
    \deligne{\ZZ}{1} \qi \sho{X}^\unit [-1]
  \end{math}
  via the standard exponential sequence, so that
  \begin{math}
    \delH[k]{X}{\ZZ}{1} \iso H^{k -1}(X, \sho{X}^\unit)\,.
  \end{math}
  In particular
  \begin{math}
    \delH[1]{X}{\ZZ}{1} \iso H^0(X,\sho{X}^\unit)\,,
  \end{math}
  the global invertibles on $X$, and
  \begin{math}
    \delH[2]{X}{\ZZ}{1} \iso \pic (X)\,,
  \end{math}
  the Picard group of line bundles over $X$.
\end{example}
\begin{example}
  \begin{math}
    \deligne{\ZZ}{2} \qi \bigl( \sho{X}^\unit
    \xrightarrow{d\log} \shomega[1]{X}\bigr)[-1]
  \end{math}\,.
  A fundamental observation by Deligne (see
  ref.~\cite{MR86h:11103}) is that $\delH[2]{X}{\ZZ}{2}$ is
  identified with the group of isomorphism classes of holomorphic
  line bundles with (holomorphic) connection. This is easily
  understood from a \cech\ cohomology point of view. Using the
  cover $\cover{U}_X = \{U_i\}_{i\in I}$, a class in
  \begin{displaymath}
    \delH[2]{X}{\ZZ}{2} \iso
    \HHH^1(X,\sho{X}^\unit
    \xrightarrow{d\log} \shomega[1]{X})
  \end{displaymath}
  is represented by a pair
  \begin{math}
    (\omega_i, g_{ij}) 
  \end{math}
  with 
  \begin{math}
    \omega_i \in \shomega[1]{X}(U_i)
  \end{math}
  and 
  \begin{math}
    g_{ij} \in \sho{X}^\unit (U_{ij})
  \end{math}
  satisfying the relations
  \begin{displaymath}
    \omega_j - \omega_i = d\log g_{ij}\,,\quad
    g_{ij} g_{jk} = g_{ik}\,.
  \end{displaymath}
  The \cech\ representative for the actual class in
  $\delH[2]{X}{\ZZ}{2}$ is obtained (up to a multiplication by
  $\tate$) by extracting local logarithms $\log g_{ij}$, see
  ref.~\cite{esn-vie:del} for full details.
\end{example}

For \emph{real} Deligne cohomology, i.e. when $A=\RR$, other
models quasi-isomorphic to those in eqs.~\eqref{eq:1} and
\eqref{eq:2} are available. Since the maps
\begin{displaymath}
  \bigl( \RR(j) \rightarrow \shomega{X}\bigr)
  \lqi \bigl( \RR(j) \rightarrow \CC \bigr)
  \lqi \RR(j-1)
  \lqi \she{X}(j-1)
\end{displaymath}
are all quasi-isomorphisms in the derived category, cf.
\cite{esn-vie:del}, we have
\begin{equation}
  \label{eq:5}
  \deligne{\RR}{j}
  \lqi
  \cone \big(F^j\shomega{X} \rightarrow \she{X}(j-1)\big)[-1]\,.
\end{equation}

Moreover, we can use smooth forms thanks to the fact that the
inclusion $\shomega{X} \hookrightarrow \sha{X}$ is a filtered
quasi-isomorphism with respect to the filtrations $F^j\shomega{X}
\hookrightarrow F^j\sha{X}$. Here $F^j\sha{X}$ is the subcomplex
of $\sha{X}$ comprising forms of type $(p,q)$ where $p$ is at
least $j$, so that $F^j\sha[n]{X} = \oplus_{p\geq j}
\sha[{p,n-p}]{X}$.

Let $(\omega\sb 1, \eta\sb 1)$ be an element of degree $n$ in
$\deligne{\RR}{j}$---this means that $\omega\sb 1\in
F^j\shomega[n]{X}$ and $\eta\sb 1\in \she[n-1]{X}(j-1)$---and
$(\omega\sb 2, \eta\sb 2)$ any element in $\deligne{\RR}{k}$.  A
product is given by the formula (cf.  ref.~\cite{esn-vie:del}):
\begin{equation}
  \label{eq:6}
  (\omega_1,\eta_1) \,\Tilde\cup\, (\omega_2,\eta_2)
  = \bigl(\omega\sb 1 \wedge \omega\sb 2 ,
  (-1)^n\, \pi\sb p\omega\sb 1 \wedge \eta\sb 2
  +\eta\sb 1\wedge \pi\sb q \omega\sb 2 \bigr)\,.
\end{equation}
\begin{example}
  $\delH[1]{X}{\RR}{1}$ is the group of real valued functions
  $\eta$ on $X$ such that there exists a holomorphic one-form
  $\omega$ such that $\pi\sb 0\omega = d\eta$. In other words, it
  is the group of those real smooth functions $\eta$ such that
  $\del\eta$ is holomorphic. In particular, if $f$ is holomorphic
  and invertible on $U \subset X$, then the class in
  $\delH[1]{X}{\RR}{1}$ determined by $f$ is represented by
  $(d\log f, \log \abs{f})$.
\end{example}

\subsection{Cones}
\label{sec:cones}

We recall here a variant of \bei's formula for the cup product on
certain diagrams of complexes. (For full details see
refs.~\cite{math.CV/0211055,bei:hodge_coho,esn-vie:del}.)

For $i=1,2,3$ consider the diagrams of complexes
\begin{equation}
  \label{eq:7}
  \mathcal{D}_i \eqdef
  X^\bullet_i \overset{f_i}{\lto} Z^\bullet_i
  \overset{g_i}{\longleftarrow} Y^\bullet_i
\end{equation}
and set
\begin{displaymath}
  C(\mathcal{D}_i) =
  \cone (X^\bullet_i\oplus Y^\bullet_i
  \xrightarrow{f_i-g_i} Z^\bullet_i)[-1]\,,
  \quad i=1,2,3\,.
\end{displaymath}
Suppose there are product maps $X^\bullet_1\otimes X^\bullet_2
\xrightarrow{\cup} X^\bullet_3$, and similarly for $Y^\bullet_i$,
and $Z^\bullet_i$. We assume the products to be compatible with
the $f_i$, $g_i$ only up to homotopy, namely there exist maps
\begin{displaymath}
  h \colon \bigl( X_1\otimes X_2 \bigr)^\bullet
  \lto Z_3^{\bullet -1}\:,\quad
  k \colon \bigl( Y_1\otimes Y_2 \bigr)^\bullet
  \lto Z_3^{\bullet -1}
\end{displaymath}
such that
\begin{displaymath}
  f_3\circ \cup - \cup \circ (f_1\otimes f_2)
  = d\, h + h\,d\:, \quad
  g_3\circ \cup - \cup \circ (g_1\otimes g_2)
  = d\, k + k\,d\,,
\end{displaymath}
with obvious meaning of the symbols. The following lemma
establishes a variant of \bei's product
formula~\cite{bei:hodge_coho}.
\begin{lemma}
  For
  \begin{math}
    (x_i,y_i,z_i) \in X^\bullet_i \oplus Y^\bullet_i \oplus
    Z^{\bullet -1}_i\,,\;i=1,2\,,
  \end{math}
  and a real parameter $\alpha$, the following formula:
  \begin{equation}
    \label{eq:8}
    \begin{split}
      (x_1,y_1,z_1) \cup_\alpha (x_2,y_2,z_2) =
      \Big(&x_1\cup x_2, y_1\cup y_2, \\
      &(-1)^{\deg (x_1)}
      \big((1-\alpha )f_1(x_1) + \alpha g_1(y_1) \big) \cup z_2 \\
      &\quad +z_1\cup \big( \alpha f_2(x_2) +
      (1-\alpha)g_2(y_2)\big)\\
      &\qquad -h(x_1\otimes x_2) +k(y_1\otimes y_2)
      \Big)\,.
    \end{split}   
  \end{equation}
  defines a family of products
  \begin{displaymath}
    C(\mathcal{D}_1)\otimes C(\mathcal{D}_2)
    \xrightarrow{\cup_\alpha}
    C(\mathcal{D}_3)\,.
  \end{displaymath}
  These products are homotopic to one another, and graded
  commutative up to homotopy. The homotopy formula is the same as
  that found in ref.~\cite{bei:hodge_coho}.
\end{lemma}
\begin{proof}
  Direct verification.
\end{proof}
If the maps $f_i$, $g_i$ above are strictly compatible with the
products, namely the homotopies $h$ and $k$ are zero,
\eqref{eq:8} reduces to the formulas found in
\cite{bei:hodge_coho,esn-vie:del}. Homotopy commutativity at the
level of complexes ensures the corresponding cohomologies will
have genuine graded commutative products.

\section{Hermitian holomorphic Deligne cohomology}
\label{sec:herm-holom-deligne}

\subsection{Metrized line bundles}
\label{sec:metr-line-bundl}

Let $X$ be a complex manifold. Consider a holomorphic line bundle
$L$ on $X$ with hermitian fiber metric $\rho$ or, equivalently,
an invertible sheaf $L$ equipped with a map
\begin{math}
  \rho \colon L \to \she[0]{X,+}
\end{math}
to (the sheaf of) positive real smooth functions, see
ref.~\cite{lang:arakelov} for the relevant formalism. Let
$\widehat{ \pic (X)}$ denote the group of isomorphism classes of
line bundles with hermitian metric. A basic observation by
Deligne (cf.~\cite{esn:char}) is that $\widehat{\pic X}$ can be
identified with the second hypercohomology group:
\begin{equation}
  \label{eq:9}
  \HHH^2\bigl(X,%
  \ZZ(1)_X \overset{\imath}{\lto} \sho{X}
  \xrightarrow{-\pi_0} \she[0]{X}\bigr)\,.
\end{equation}
This is easy to see in \cech\ cohomology. Suppose $s_i$ is a
trivialization of $L\rvert_{U_i}$, with transition functions
$g_{ij}\in \sho{X}^\unit (U_{ij})$ determined by $s_j = s_i
g_{ij}$. Let $\rho_i$ be the value of the quadratic form
associated to $\rho$ on $s_i$, namely $\rho_i = \rho(s_i)$. Then
we have $\rho_j = \rho_i \abs{g_{ij}}^2$. Taking logarithms, we
see that
\begin{equation*}
  \bigl(\tate c_{ijk},
  \log g_{ij}, \tfrac{1}{2} \log \rho_i \bigr) \,, 
\end{equation*}
where
\begin{math}
  \tate c_{ijk} = \log g_{jk} -\log g_{ik} +\log g_{ij}\in
  \ZZ(1)\,,
\end{math}
is a cocycle representing the class of the pair $(L,\rho)$.

\subsubsection{Canonical connection}
\label{sec:canonical-connection}

Recall for later use that the \emph{canonical
  connection,}~\cite{gh:alg_geom} on a metrized line bundle
$(L,\rho)$ is the unique connection compatible with both the
holomorphic and hermitian structures.  In \cech\ cohomology with
respect to the cover $\cover{U}_X$ as above, the canonical
connection on $(L,\rho)$ corresponds to a collection of $(1,0)$
forms
\begin{math}
  \xi_i\in \sha[{1,0}]{X}(U_i)
\end{math}
satisfying the relations
\begin{align}
  \label{eq:10} \xi_j - \xi_i &= d\log g_{ij}\\
  \label{eq:11} \pi_0(\xi_i) &= \tfrac{1}{2} d\log \rho_i\,.
\end{align}
The latter just means 
\begin{math}
  \xi_i = \del \log \rho_i\,,
\end{math}
in more familiar terms. The global $2$-form
\begin{equation}
  \label{eq:12}
  c_1(\rho) = \eta_i \coin \delb\del\log \rho_i
\end{equation}
represents the first Chern class of $L$ in $H^2(X,\RR(1))$. The
class of $c_1(\rho)$ is in fact a pure Hodge class in
$H^{1,1}(X)$---the image of the first Chern class of $L$ under
the map $\delH[2]{X}{\ZZ}{1} \to H^2_\mathcal{D}(X,\RR(1))$
induced by $\ZZ(1) \to \RR(1)$. It only depends on the class of
$(L,\rho)$ in $\widehat{\pic (X)}$.

\subsection{Hermitian holomorphic complexes}
\label{sec:herm-holom-compl}

In ref.~\cite{bry:quillen} Brylinski introduced the complexes
\begin{equation}
  \label{eq:13}
  C(l)^\bullet = \cone \bigl(
  \ZZ(l)_X \oplus (F^l\!\sha{X}\cap \sigma^{2l}\she{X}(l))
  \lto \she{X}(l)
  \bigr)[-1]\,.
\end{equation}
\begin{definition}
  The hypercohomology groups
  \begin{equation}
    \label{eq:14}
    \dhhH[p]{X}{l} \eqdef \HHH^p(X,C(l))
  \end{equation}
  are the \emph{Hermitian holomorphic Deligne} cohomology groups.
\end{definition}
By the remark after eq.~\eqref{eq:5}, the complex
\begin{displaymath}
  \deltilde{\RR}{l} =
  \cone \big(F^l\!\sha{X} \rightarrow \she{X}(l-1)\big)[-1]\,.
\end{displaymath}
also computes the real Deligne cohomology. Then consider the
complex
\begin{equation}
  \label{eq:15}
  \dhh{l} = \cone \bigl(
  \deligne{\ZZ}{l}\oplus (F^l\!\sha{X}\cap \sigma^{2l}\she{X}(l))
  \lto \deltilde{\RR}{l}
  \bigr)[-1]\,.
\end{equation}
In ref.~\cite{math.CV/0211055} we prove
\begin{lemma}
  The complexes $C(l)^\bullet$ and $\dhh{l}$ are
  quasi-isomorphic, hence we also have
  \begin{displaymath}
    \dhhH[p]{X}{l} = \HHH^{\,p}(X,\dhh{l})\,.
  \end{displaymath}
\end{lemma}
\begin{remark}
  The complex 
  \begin{math}
    F^l\!\sha{X}\cap \sigma^{2l}\she{X}(l)
  \end{math}
  appearing in both~\eqref{eq:13} and \eqref{eq:14} can be
  rewritten in terms of the complex $G(l)^\bullet$ of
  ref.~\cite{esn:char}. Set
  \begin{displaymath}
    G(l)^\bullet = 0\lto \dotsm \lto 0
    \lto \sha[{(l,l)}]{X}
    \overset{d}{\lto}
    \sha[{(l+1,l)}]{X}\oplus \sha[{(l,l+1)}]{X}
    \overset{d}{\lto}\dotsm \;.
  \end{displaymath}
  Then we have
  \begin{math}
    F^l\!\sha{X}\cap \sigma^{2l}\she{X}(l)
    = G(l)^\bullet \cap \she{X}(l)\,.
  \end{math}
\end{remark}
For certain ranges of values of the cohomology index the groups
$\dhhH[p]{X}{l}$ are fairly ordinary. Indeed we have the
following easy
\begin{lemma}
  \label{lem:2}
  For $p\leq 2l-1$ we have
  \begin{displaymath}
    \dhhH[p]{X}{l} \iso H^{p-1}(X,\RR(l)/\ZZ(l))\,.
  \end{displaymath}
\end{lemma}
\begin{proof}
  Using either $C(l)^\bullet$ or $\dhh{l}$, we see that they are
  quasi-isomorphic to 
  \begin{displaymath}
    \cone \bigl(
    F^l\!\sha{X}\cap \sigma^{2l}\she{X}(l) \lto
    \RR(l)/\ZZ(l)
    \bigr)[-1]\,,
  \end{displaymath}
  which leads to the triangle
  \begin{displaymath}
    \RR(l)/\ZZ(l)[-1] \lto \dhh{l} \lto
    F^l\!\sha{X}\cap \sigma^{2l}\she{X}(l)
    \overset{+1}{\lto}\,.
  \end{displaymath}
  The statement follows.
\end{proof}
In general these groups are interesting when $p\geq 2l$. The most
important example is:
\begin{lemma}
  \begin{displaymath}
    \widehat{\pic (X)} \iso \dhhH[2]{X}{1}\,.
  \end{displaymath}
\end{lemma}
\begin{proof}
  We have quasi-isomorphisms
  \begin{displaymath}
    \ZZ(1)_X \overset{\imath}{\lto} \sho{X}
    \xrightarrow{-\pi_0} \she[0]{X}
    \lqi \dhh{1} \lqi C(1)^\bullet\,.
  \end{displaymath}
  Indeed, note that $\dhh{1}$ can be rewritten as
  \begin{displaymath}
    \cone \bigl( \deligne{\ZZ}{1}
    \to \deltilde{\RR}{1}/ (F^1\!\sha{X}\cap
    \sigma^{2}\she{X}(1))\bigr)[-1]
  \end{displaymath}
  and
  \begin{displaymath}
    \deltilde{\RR}{1}/ (F^1\!\sha{X}\cap
    \sigma^{2}\she{X}(1))
    \lqi
    \cone \bigl(
    F^1\! \sha{X} / F^1\!\sha{X}\cap \sigma^{2}\she{X}(1)
    \xrightarrow{-\pi_0} \she{X}
    \bigr)[-1]\,.
  \end{displaymath}
  By direct verification, the latter complex is quasi-isomorphic
  to $\she[0]{X}[-1]$. Thus
  \begin{displaymath}
    \dhh{1} \lqi
    \cone \bigl( \deligne{\ZZ}{1} \to \she[0]{X}[-1] \bigr)[-1]
    \lqi \ZZ(1)_X \to \sho{X} \to \she[0]{X}\,.
  \end{displaymath}
\end{proof}
Since hermitian holomorphic Deligne complexes can be expressed as
cones of diagrams of the form \eqref{eq:7}, they admit cup
products, and hence there is a cup product for hermitian
holomorphic Deligne cohomology \cite{bry:quillen}:
\begin{displaymath}
  \dhhH[p]{X}{l} \otimes \dhhH[q]{X}{k}
  \overset{\cup}{\lto}
  \dhhH[p+q]{X}{l+k}\,.
\end{displaymath}

\subsection{Explicit cocycles}
\label{sec:explict-cocycles}

Use of the seemingly more complicated complex \eqref{eq:15} in
place of the one in~\eqref{eq:13} is justified by the fact that
the data comprising the canonical connection can be characterized
cohomologically, as follows:
\begin{lemma}
  Let $(L,\rho)$ be a metrized line bundle on $X$. Assume
  $(L,\rho)$ to be trivialized with respect to the open cover
  $\cover{U}_X$ of $X$ as before. The data:
  \begin{gather*}
    \xi_i\in \sha[(1,0)]{X}(U_i)\,,\quad
    \tfrac{1}{2} \log\rho_i \in \she[0]{X}(U_i)\,,\quad
    \eta_i\in \sha[(1,1)]{X}(U_i)\,,\\
    \tate c_{ijk}\in \ZZ(1)_X (U_{ijk})\,, \quad
    \log g_{ij}\in \sho{X}(U_{ij})
  \end{gather*}
  represent a degree $2$ cocycle with values in 
  \begin{math}
    \Tot \vC{\cover{U}_X}{\dhh{1}}
  \end{math}
  if and only if the relations \eqref{eq:10}, \eqref{eq:11},
  \eqref{eq:12}, plus those in sect.~\ref{sec:metr-line-bundl},
  defining the canonical connection are satisfied.
\end{lemma}
\begin{proof}
  One need only unravel the cone defining $\dhh{1}$ as follows:
  \begin{equation}
    \label{eq:16}
    \begin{CD}
      \ZZ(1)_X @>>> \sho{X} @>>>  0 @>>> \dotsm \\
      & & @VV{0\oplus\pi_0}V @VVV    \\
      & & F^1\!\sha[1]{X} \oplus \she[0]{X} @>>> F^1\!\sha[2]{X}
      \oplus \she[1]{X} @>>> \dotsm \\
      & & & & @AA{\jmath\oplus 0}A & \\
      & & & & F^1\!\sha[2]{X}\cap \she[2]{X}(1) @>>> \dotsm
    \end{CD}
  \end{equation}
  and carefully chase the diagram.
\end{proof}
On the other hand, the hermitian holomorphic Deligne complex in
the form \eqref{eq:13} corresponds to ``reducing the structure
group'' from $\CC^\unit$ to $\TT$. This can be made explicit for
$l=1$ and a line bundle $L\to X$ by choosing sections $t_i$ of
the smooth bundle corresponding to $L$ such that $\rho(t_i)=1$.
Clearly the resulting smooth transition functions will be
sections of $\sheaf{\TT}_X$ over $U_{ij}$. See
refs.~\cite{bry:quillen} and \cite{brymcl:deg4_II} for more
details.

\section{Tame symbol and hermitian structure}
\label{sec:tame-symb-herm}

Let $X$ be a complex analytic manifold and $U\subset X$ open. Let
$f$ and $g$ two invertible holomorphic functions on $U$. The 
tame symbol~\cite{del:symbole} $\tame{f}{g}$ associated to $f$ and
$g$ is a $\sho{X}^\unit\vert_U$-torsor equipped with an analytic
connection. 

\subsection{Cup product and Deligne torsor}
\label{sec:cup-product}

(See \cite{del:symbole,esn-vie:del}.) We consider $f$ and $g$ as
elements of $\delH[1]{U}{\ZZ}{1}$. Then
\begin{math}
  \tame{f}{g} = f\cup g \in \delH[2]{U}{\ZZ}{2}\,.
\end{math}
Consider the cover $\cover{U}_X$ of $X$ so that $U$ is covered by
$\{U\cap U_i\}_{i\in I}$ and choose representatives 
\begin{math}
  (\tate\,m_{ij}, \log_i\! f)
\end{math}
and
\begin{math}
  (\tate\,n_{ij}, \log_i\! g)
\end{math}
for $f$ and $g$, respectively. Then, using \eqref{eq:4}, the cup
product is represented by the cocycle:
\begin{equation}
  \label{eq:17}
  \Bigl( (\tate)^2 m_{ij}n_{jk}\,,\:
  -\tate\,m_{ij}\, \log_j\! g\,,\:
  \log_i\! f\: \frac{d g}{g}\,\Bigr)\,.
\end{equation}
Under the quasi-isomorphism with the complex
\begin{math}
  \bigl(\sho{X}^\unit \to \shomega[1]{X}\bigr)
\end{math}
(which essentially amounts to a division by $\tate$)
the cocycle~\eqref{eq:17} becomes
\begin{equation}
  \label{eq:18}
  \big( g^{-m_{ij}}\,,
  -\frac{1}{\tate}\log_i\! f\: \frac{dg}{g} \bigr)\,.
\end{equation}
In ref.~\cite{del:symbole} the trivializing section on $U\cap
U_i$ corresponding to~\eqref{eq:18} is denoted
\begin{math}
  \lbrace \log_i\! f,g\rbrace\,.
\end{math}
Two trivializations over $U\cap U_i$ and $U\cap U_j$ are related
by 
\begin{math}
  \lbrace \log_j\! f,g\rbrace
  = \lbrace \log_i\! f,g\rbrace\,
  g^{-m_{ij}}\,.
\end{math}
Furthermore, the analytic connection is defined by the rule:
\begin{equation}
  \label{eq:19}
  \nabla \lbrace \log_i\! f,g\rbrace = -
  \lbrace \log_i\! f,g\rbrace \otimes
  \frac{1}{\tate}\log_i\! f\: \frac{dg}{g}\,.
\end{equation}
A general section $s$ of $\tame{f}{g}$ can be written as 
\begin{math}
  s = h_i \, \lbrace \log_i\! f,g\rbrace \,,
\end{math}
for some 
\begin{math}
  h_i \in{ \sho{U}(U_i)}\,,
\end{math}
and therefore
\begin{equation}
  \label{eq:20}
  \nabla s = \lbrace \log_i\! f,g\rbrace \otimes \bigl( dh_i
  -\frac{1}{\tate}\log_i\! f\: \frac{dg}{g} \bigr)\,.
\end{equation}

\subsection{Heisenberg group}
\label{sec:heisenberg-group}

An equivalent approach to the Deligne symbol is via the complex
three-dimensional Heisenberg group, see
refs.~\cite{bloch:dilog_lie,MR94k:19002,rama:reg_hei}. Let
$H_\CC$ denote the group of complex unipotent $3\times 3$ lower
triangular matrices. Let
\begin{displaymath}
  H_\ZZ = \left\lbrace
    \begin{pmatrix}
      1&&\\
      m_1&1&\\
      m_2&n_1&1
    \end{pmatrix}
    \Big\vert\; m_1,n_1\in \ZZ(1)\,,\; m_2\in\ZZ(2) \right\rbrace
  \subset H_\CC\,.
\end{displaymath}
The quotient $H_\CC / H_\ZZ$ is a $\CC/\ZZ(2)$-bundle over 
\begin{math}
  \CC / \ZZ(1)\times \CC / \ZZ(1)
\end{math}
via the projection map
\begin{displaymath}
  p\colon
  \begin{bmatrix}
    1&&\\
    x&1&\\
    z&y&1
  \end{bmatrix} \mapsto ([x],[y])\,,
\end{displaymath}
where $x,y,z\in \CC$, and the brackets denote the appropriate
equivalence classes. (The $\CC/\ZZ(2)$-action is by
multiplication with a matrix of the form
\begin{math}
  \Bigl( 
  \begin{smallmatrix}
    1&&\\ 0&1&\\z&0&1
  \end{smallmatrix}
  \Bigr)\,.
\end{math})

The twisting of $H_\CC/H_\ZZ$ is analogous to that of the Deligne
torsor in sect.~\ref{sec:cup-product}: the right action of
$H_\ZZ$ on $H_\CC$ amounts to:
\begin{equation}
  \label{eq:21}
  x \mapsto x + m_1\,,\quad
  y \mapsto y + n_1\,,\quad
  z \mapsto z + m_1\cdot y + m_2\,.
\end{equation}
Moreover, the complex form
\begin{equation}
  \label{eq:22}
  \omega = \frac{1}{\tate}(dz - x\,dy)
\end{equation}
is invariant under the action of $H_\ZZ$ and defines a
$\CC/\ZZ(2)$-connection form on the total space $H_\CC/H_\ZZ$.

The invertible functions $f$ and $g$ on $U$ define a map
$(f,g)\colon U \to \CC^\unit \times \CC^\unit$. Then
the tame symbol $\tame{f}{g}$ is obtained as the pull-back:
\begin{displaymath}
  \tame{f}{g}= (f,g)^* \bigl( H_\CC/H_\ZZ\bigr)\,,
\end{displaymath}
and the section $\lbrace\log_i\!f,g\rbrace$ corresponds to the
class of the matrix
\begin{displaymath}
  \begin{pmatrix}
    1 & &\\
    \log_i\!f & 1 & \\
    0 & \log_i\!g & 1
  \end{pmatrix} \,.
\end{displaymath}
Furthermore, the pull-back of the connection form $\omega$ on
$H_\CC/H_\ZZ$ along the section $\lbrace \log_i\!f,g\rbrace$ is
the same form as the one in~\eqref{eq:17}. More generally, a
section $s$ as given at the end of sect.~\ref{sec:cup-product}
corresponds to the class of the matrix
\begin{displaymath}
  \begin{pmatrix}
    1 & &\\
    \log_i\!f & 1 & \\
    h_i & \log_i\!g & 1
  \end{pmatrix} \,,
\end{displaymath}
Pulling back~\eqref{eq:22} along the section gives~\eqref{eq:20}.

\subsection{Hermitian product structure}
\label{sec:hermitian-structure}

Consider the ``imaginary part'' map
\begin{equation}
  \label{eq:23}
  \begin{aligned}
    \CC \otimes \CC &\lto \RR(1)\\
    a \otimes b &\longmapsto -\pi_1(a)\,\pi_0(b)
    \coin -\sqrt{-1} \Im (a) \Re (b)\,,
  \end{aligned}
\end{equation}
Similarly, we have:
\begin{equation}
  \label{eq:24}
    \sho{X}\otimes \sho{X} \lto \she[0]{X}(1)\,\quad
    f \otimes g \longmapsto -\pi_1(f)\,\pi_0(g)\,.
\end{equation}
\begin{definition}
  \label{def:1}
  Define the map
  \begin{equation}
    \label{eq:25}
    \begin{split}
      \bigl( \ZZ(1)_X \to \sho{X} \bigr)\otimes
      \bigl( \ZZ(1)_X \to \sho{X} \bigr)
      &\lto
      \bigl( \ZZ(2)_X \to \sho{X} \xrightarrow{-\pi_1}
      \she[0]{X}(1) \bigr)\\
      &\lqi \tate\otimes
      \bigl( \ZZ(1)_X \to \sho{X} \xrightarrow{-\pi_0}
      \she[0]{X} \bigr)
    \end{split}
  \end{equation}
  by using~\eqref{eq:24} in place of the map
  \begin{math}
    \sho{X}\otimes \sho{X} \to \shomega[1]{X}\,,
  \end{math}
  \begin{math}
    f\otimes g \mapsto fdg\,,
  \end{math}
  in \eqref{eq:4}.
\end{definition}
\begin{proposition}
  \label{prop:1}
  The product map~\eqref{eq:25} is well defined, namely it is a
  map of complexes. Furthermore, it is homotopy graded
  commutative.
\end{proposition}
\begin{proof}
  The fact that~\eqref{eq:25} is a map of complexes is a direct
  verification. After ref.~\cite{esn-vie:del}, consider the map
  \begin{equation*}
    h (f\otimes g) = f\,g\,,\quad f,g\in\sho{X}\,,
  \end{equation*}
  and zero otherwise. It provides the required homotopy.
\end{proof}
The target complex of the product map in eq.~\eqref{eq:25} is the
complex encoding hermitian structures appearing in
sect.~\ref{sec:metr-line-bundl}. In other words, up to
quasi-isomorphism, we have a product:
\begin{equation*}
  \deligne{\ZZ}{1} \otimes \deligne{\ZZ}{1}
  \lto \tate\otimes\dhh{1}\,.
\end{equation*}
\begin{remark}
  The map~\eqref{eq:24} provides an explicit homotopy map for the
  homotopy commutative diagram
  \begin{equation*}
    \begin{CD}
      \deligne{\ZZ}{1}\otimes \deligne{\ZZ}{1} @>>>
      \deligne{\ZZ}{2}\\
      @VVV @VVV\\
      \deligne{\RR}{1}\otimes \deligne{\RR}{1} @>>>
      \deligne{\RR}{2}
    \end{CD}
  \end{equation*}
  where the model~\eqref{eq:5} for $\deligne{\RR}{k}$ is used
  (see \cite{esn-vie:del}).
\end{remark}
Now, in view of Prop.~\ref{prop:1}, we have a graded commutative
product at the level of cohomology groups. In particular, let
$f,g$ be two invertible holomorphic functions on $U\subset X$.
\begin{proposition}
  \label{prop:2}
  The Deligne torsor underlying $\tame{f}{g}$ admits a hermitian
  fiber metric.
\end{proposition}
\begin{proof}
  View $f$ and $g$ as elements of $\delH[1]{U}{\ZZ}{1}$. Taking
  the product according to~\eqref{eq:25} yields an element in
  \begin{equation*}
    \dhhH[2]{U}{1}\iso \widehat{\pic (U)}
  \end{equation*}
  that is, a holomorphic line bundle with hermitian fiber metric
  (up to isomorphism).
  
  Taking the image of the tame symbol $\tame{f}{g}$ under the map
  \begin{math}
    \delH{U}{\ZZ}{2} \to \delH{U}{\ZZ}{1}=\pic (U)
  \end{math}
  induced by
  \begin{math}
    \deligne{\ZZ}{2} \to \deligne{\ZZ}{1}
  \end{math}
  forgets the analytic connection and retains just the line
  bundle. Similarly, the map
  \begin{math}
    \dhhH[2]{U}{1} \to \delH{U}{\ZZ}{1}=\pic (U)
  \end{math}
  induced by
  \begin{math}
    \dhh{1} \to \deligne{\ZZ}{1}
  \end{math}
  forgets the hermitian structure. Clearly both map to the same
  underlying line bundle.
\end{proof}
Using a \cech\ cover we can represent $f$ and $g$ as in
sect.~\ref{sec:cup-product}. Then the cocycle corresponding to
their product in
\begin{math}
  \dhhH[2]{U}{1}
\end{math}
is:
\begin{equation}
  \label{eq:26}
  \Bigl( \tate\, m_{ij}\,n_{jk}\,,\:
  -m_{ij} \log_j\! g\,,
  -\frac{1}{\tate}\,\pi_1(\log_i\! f) \log\abs{g}\,\Bigr)\,.
\end{equation}
This allows us to identify the representative of the hermitian
metric, or rather its logarithm, as
\begin{equation}
  \label{eq:27}
  \onehalf \log\rho_i =
  -\frac{1}{\tate}\,\pi_1(\log_i\! f) \log\abs{g}\,.
\end{equation}
It follows that if $s$ is the local section at the end of sect.\
\ref{sec:cup-product} then
\begin{equation}
  \label{eq:67}
  \log \rho (s) = \frac{1}{\tate}
  \bigl(\pi_1(h_i)-\pi_1(\log_i\!f)\,\log \abs{g}\bigr)\,.
\end{equation}

\subsubsection{Remarks on the Heisenberg bundle}
\label{sec:remarks-heis-group}

The hermitian metric can be constructed from the more global
point of view afforded by the use of the Heisenberg group
recalled in sect.~\ref{sec:heisenberg-group}. The hermitian
metric on the bundle
\begin{math}
  H_\CC /H_\ZZ \to \CC^\unit \times \CC^\unit
\end{math}
is given by the map $\rho : H_\CC/H_\ZZ\to \RR_+$ defined by:
\begin{equation}
  \label{eq:28}
  \rho\colon
  \begin{bmatrix}
    1&&\\
    x&1&\\
    z&y&1
  \end{bmatrix} 
  \longmapsto
  \exp \frac{1}{\tate}
  \bigl(\pi_1(z)-\pi_1(x)\,\pi_0(y)\bigr)
\end{equation}
Indeed, using the explicit
action~\eqref{eq:21}, one checks~\eqref{eq:28} is invariant and
provides the required quadratic form. In particular, the
quantity
\begin{equation*}
  -\frac{1}{\tate}\,\pi_1(x)\,\pi_0(y)
\end{equation*}
is immediately shown to behave as the logarithm of the local
representative of a hermitian metric. Thus the hermitian
holomorphic line bundle represented by the cocycle~\eqref{eq:26}
is the pull-back of $(H_\CC/H_\ZZ,\rho)$ via the map $(f,g)\colon
U \to \CC^\unit \times \CC^\unit$.

\subsubsection{Relations with Mixed Hodge Structures}
\label{sec:relations-with-mixed}

Both structures, namely the standard cup product
\begin{math}
  \deligne{\ZZ}{1} \otimes \deligne{\ZZ}{1}\to \deligne{\ZZ}{2}
\end{math}
given by~\eqref{eq:4}, and the modified one
\begin{math}
  \deligne{\ZZ}{1} \otimes \deligne{\ZZ}{1}\to
  \tate\otimes \dhh{1}
\end{math}
of Definition~\ref{def:1}, can be obtained by taking projections
of a common object in two different ways.

Let $s$ be a local section of the pull-back
\begin{equation*}
  \tame{f}{g}= (f,g)^* \bigl( H_\CC/H_\ZZ\bigr)
\end{equation*}
as at the end of sect.\ \ref{sec:cup-product}.  (The local
expression in terms of matrices is given at the end of sect.\
\ref{sec:heisenberg-group}.)  Equivalently, $s$ can be considered
as a (local) lift of the map $(f,g): X \to \CC^\unit\times
\CC^\unit$ to $H_\CC/H_\ZZ$.

Let $\hodge{M}^{(2)}_X$ be the resulting variation of Mixed Hodge
Structures on $X$ obtained by pulling back the universal MHS
$\hodge{M}^{(2)}$ on $H_\CC/H_\ZZ$ via $s$.
\begin{lemma}[See \cite{MR99i:19004}]
  The period
  \begin{math}
    P(\hodge{M}^{(2)}_X) \in \sho{X}\otimes_\QQ \sho{X}
  \end{math}
  of $\hodge{M}^{(2)}_X$ is given by:
  \begin{equation*}
    \begin{split}
      P(\hM^{(2)}_X) &= \frac{h}{(\tate)^2}\otimes 1
      -1\otimes \frac{h}{(\tate)^2}\\
      &+1\otimes \frac{\log f\,\log g}{(\tate)^2}
      -\frac{\log g}{\tate} \otimes \frac{\log f}{\tate}
    \end{split}
  \end{equation*}
\end{lemma}
\begin{proof}
  The expression is computed in the appendix for the universal case.
\end{proof}
Notice that the period actually belongs to the kernel of the
multiplication map $a\otimes b\to ab$.

Let us now use the map
\begin{math}
  \sho{X}\otimes_\QQ \sho{X} \to \sho{X}\otimes_\CC \sho{X}\,.
\end{math}
Let $\curly{I}_X$ be the kernel of the multiplication map (over
\CC). Then $\shomega[1]{X/\CC}\iso
\curly{I}_X/\curly{I}_X^2$. The calculations for the following
proposition are done in the universal case in the appendix.
\begin{proposition}
  The expressions~\eqref{eq:20} and~\eqref{eq:67} respectively
  correspond to the images of $P(\hM^{(2)}_X)$ under the projections
  \begin{math}
    \curly{I}_X\subset \sho{X}\otimes_\CC \sho{X} \to
    \shomega[1]{X/\CC}\,,
  \end{math}
  sending $a\otimes b - ab\otimes 1$ to $a\,db$, and
  \begin{math}
    \curly{I}_X\subset \sho{X}\otimes_\CC \sho{X} \to \she[0]{X}
  \end{math}
  given by~\eqref{eq:24}.
\end{proposition}

\subsection{Comparisons}
\label{sec:comparisons}

In the previous sections we have shown that the Deligne torsor
$\tame{f}{g}$ associated to two invertible functions $f$ and $g$
naturally acquires two structures: the analytic connection
$\nabla$ described in section~\ref{sec:cup-product} via the
standard cup product in Deligne cohomology, and the hermitian
structure described in section~\ref{sec:hermitian-structure} via
the modified cup product~\eqref{eq:25}. We wish to briefly
compare the two structures.

First, observe that using the canonical connection (cf.\
section~\ref{sec:canonical-connection}) a pair $(L,\rho)$ can
also be thought of as a triple $(L,\rho, \nabla^\rho)$, where
$\nabla^\rho$ is the canonical connection determined by
$\rho$. Equivalently, we can just consider the pair
$(L,\nabla^\rho)$. Also, let us stress that the canonical
connection is only a \emph{smooth} connection and is in general
far from being analytic (or algebraic).

Thus our question can be reformulated as follows: for a given
line bundle $L$ equipped with an analytic connection $\nabla$ and
a hermitian fiber metric $\rho$, how do the pairs $(L,\nabla)$
and $(L,\nabla^\rho)$ compare?

The answer is the following well-known
\begin{lemma}
  \label{lem:1}
  Consider both $\nabla$ and $\nabla^h$ as \emph{smooth}
  connections. Then:
  \begin{enumerate}
  \item $\nabla - \nabla^h$ determines a global section of
    $\sha[1,0]{X}$, and
  \item this global section is zero, that is, $\nabla=\nabla^h$,
    if and only if $L$ is unitary flat, namely it defines an
    element of $H^1 (X, \RR/\ZZ)$.
  \end{enumerate}
\end{lemma}
\begin{proof}
  It is a well-known fact that the difference of two connections
  is a global one-form. Working in a local setting, let $s\in
  L\rvert_U$ be a local section, and let $\norm{s}$ be its length
  with respect to the metric. Then $\nabla s = \omega\otimes s$,
  for $\omega \in \shomega[1]{X}(U)$, whereas $\nabla^\rho s =
  \del\log \norm{s} \otimes s$, and $\del\log \norm{s}$ gives a
  local $(1,0)$-form representative of $\nabla^\rho$, cf.\
  section~\ref{sec:canonical-connection}. Clearly, the difference
  $\omega - \del\log \norm{s}$ gives a global section of
  $\sha[1,0]{X}$.

  As for the second point, one would have $\delb\del\log
  \norm{s}=0$, but this represents $c_1({L})$, hence the
  conclusion.
\end{proof}
In the situation when the two connections agree, that is, the
connection is simultaneously analytic and it is the canonical
connection associated to a hermitian structure, we say they are
\emph{compatible.} The line bundle supporting it is necessarily
flat.

Interestingly enough, the previous lemma can be recast into
entirely cohomological terms. This is advantageous in dealing
with the special case $L=\tame{f}{g}$ of special interest to us,
as well as to address the very same question in the case of
gerbes later on in this paper.

In the previous lemma we have compared $\nabla$ and $\nabla^\rho$
by mapping their respective local representatives in
$\sha[1,0]{X}$. It will be more convenient to use the sheaf of
imaginary $1$-forms instead, namely consider $\pi_1 :
\shomega[1]{X}\to \she[1]{X}(1)$ and $d : \she[0]{X}(1)\to
\she[1]{X}(1)$. Consider the complex
\begin{equation*}
  \Lambda (2)^\bullet \eqdef
  \bigl( \ZZ(2) \overset{\imath}{\to} \sho{X}
  \xrightarrow{-\pi_1\circ d} \she[1]{X}(1) \bigr)\,,
\end{equation*}
and the obvious maps of complexes
\begin{equation*}
  \alpha : \deligne{\ZZ}{2} \lto
  \Lambda (2)^\bullet\quad
  \text{and}\quad
  \beta :  \tate\otimes\dhh{1} \lto \Lambda (2)^\bullet\,.
\end{equation*}
As usual, the cone:
\begin{equation*}
  \Gamma (2)^\bullet \eqdef
  \cone \bigl( \alpha -\beta\bigr)  [-1]\,, 
\end{equation*}
characterizes the elements in $\deligne{\ZZ}{2}$ and
$\tate\otimes\dhh{1}$ which agree in $\Lambda (2)^\bullet$. A
tedious but straightforward direct verification yields:
\begin{lemma}
  \label{lem:3}
  We have the quasi-isomorphism:
  \begin{equation}
    \label{eq:66}
    \Gamma (2)^\bullet \lqi
    \bigl( 
    \ZZ(2) \overset{\imath}{\to} \sho{X}
    \xrightarrow{(d,-\pi_1)}
    \shomega[1]{X}\oplus\she[0]{X}(1)
    \xrightarrow{\pi_1+d} \she[1]{X}(1)\bigr)
  \end{equation}
\end{lemma}
Dropping the last term in~\eqref{eq:66}, we obtain the truncation
\begin{equation*}
  \Tilde{\Gamma}(2)^\bullet \eqdef
  \bigl( 
  \ZZ(2) \overset{\imath}{\to} \sho{X}
  \xrightarrow{(d,-\pi_1)}
  \shomega[1]{X}\oplus\she[0]{X}(1)\bigr)\,,
\end{equation*}
which clearly characterizes the elements in $\deligne{\ZZ}{2}$
and $\tate\otimes\dhh{1}$ which agree in
$\tate\otimes\deligne{\ZZ}{1}$. (In other words,
$\Tilde{\Gamma}(2)^\bullet$ can be obtained by replacing $\Lambda
(2)^\bullet$ by $\deligne{\ZZ}{1}$ in the previous paragraphs.)
In particular, let us denote by $\pic (X,\nabla,h)$ the second
hypercohomology group $\HHH^2(X,\Tilde{\Gamma}(2)^\bullet)$,
namely the subgroup of $\delH[2]{X}{\ZZ}{2}\times \widehat{\pic
  (X)}$ of classes of pairs $(L,\nabla)$ and $(L,\rho)$ mapping
to the same element of $\pic (X)\iso \delH[2]{X}{\ZZ}{1}$. Then
lemma~\ref{lem:1} has the following reformulation:
\begin{lemma}
  \label{lem:4}
  There is an exact sequence:
  \begin{equation}
    \label{eq:64}
    0\lto H^1(X,\RR/\ZZ) \lto 
    \pic (X,\nabla,h) \lto E^1(X)(1)\,,
  \end{equation}
  where $E^1(X)(1)$ are the global sections of
  $\she[1]{X}(1)$. Thus compatible connections are necessarily
  flat.
\end{lemma}
\begin{proof}
   The complex $\Tilde{\Gamma}(2)^\bullet$  is a quotient of
   $\Gamma (2)^\bullet$, namely we have the exact sequence:
   \begin{equation*}
     0\lto \she[1]{X}(1) [-3] \lto \Gamma(2)^\bullet \lto
     \Tilde{\Gamma}(2)^\bullet \lto 0\,,
   \end{equation*}
   and from the resulting long exact cohomology sequence:
   \begin{equation*}
     0\to \HHH^2(X,\Gamma(2)^\bullet) \to
     \HHH^2 (X,\Tilde{\Gamma}(2)^\bullet) \to E^1(X)(1) \to
     \dotsm 
   \end{equation*}
   It was noted above that
   \begin{math}
     \HHH^2 (X,\Tilde{\Gamma}(2)^\bullet)\iso
     \pic (X,\nabla,h)\,,
   \end{math}
   whereas for $\Gamma (2)^\bullet$ we have
   \begin{equation*}
     \HHH^2 (X,{\Gamma}(2)^\bullet)\iso
     H^1(X,\RR (2)/ \ZZ(2))\,.
   \end{equation*}
   The latter isomorphism follows either from a direct
   computation, or noticing that $\Gamma (2)^\bullet$ is a
   quotient of $\dhh{2}$ (see eq.\ \eqref{eq:15}) and
   \begin{equation*}
     \HHH^2(X,\Gamma(2)^\bullet) \iso \dhhH[2]{X}{2}
   \end{equation*}
   and then using lemma~\ref{lem:2}.
\end{proof}

\subsubsection{Comparing $\protect \tame{f}{g}$ and $\protect \tamehh{f}{g}$}
\label{sec:comp-tamefg}

Suppose now $L$ is the Deligne torsor determined by two
invertible functions $f$ and $g$. Clearly, the symbols
$\tame{f}{g}$ and $\tamehh{f}{g}$ taken together determine an
element of $\pic (X,\nabla, h)$, since the underlying torsor in
$\pic (X)\iso \delH[2]{X}{\ZZ}{1}$ is the same. This element can
be represented by the cocycle
\begin{equation*}
  \Bigl(
  (\tate)^2 m_{ij}n_{jk}\,,\:
  -\tate\,m_{ij} \log_j\! g\,,\:
  \log_i\! f\, \frac{d g}{g}\oplus -\pi_1(\log_i\! f) \log\abs{g}
  \Bigr)
\end{equation*}
with values in $\Tilde{\Gamma}(2)^\bullet$.

Following Goncharov (\cite{MR1978709}) let us define for
any two invertibles $f$ and $g$ the $1$-form
\begin{equation}
  \label{eq:65}
  r_2(f,g) \eqdef \pi_1(d \log f) \log\abs{g}
  -\log\abs{f} \pi_1(d\log g)\,.
\end{equation}
This is clearly globally defined where $f$ and $g$ are
invertible.

We finally obtain the following comparison. 
\begin{proposition}
  \label{prop:5}
  The analytic connection in $\tame{f}{g}$ and the canonical one
  associated to the hermitian structure in $\tamehh{f}{g}$ are
  compatible if and only if $r_2(f,g)=0$ in $E^1(X)(1)$.
\end{proposition}
\begin{proof}
  Let $\omega_i = \log_i\!f \: dg/g$ and $\sigma_i =
  -\pi_1(\log_i\!  f) \log\abs{g}$. The connecting homomorphism
  from $\Tilde{\Gamma} (2)^\bullet$ to $\she[1]{X}(1)$, that is
  the last map to the right in the sequence~\eqref{eq:64},
  amounts to computing $\pi_1(\omega_i) + d\sigma_i$. A
  straightforward calculation yields
  \begin{equation*}
    \pi_1(\omega_i) + d\sigma_i = -r_2(f,g)\,.
  \end{equation*}
\end{proof}

\section{Hermitian holomorphic gerbes and $2$-gerbes}
\label{sec:herm-holorm-gerb}

\subsection{Higher tame symbols}
\label{sec:higher-tame-symbols}
Brylinski and McLaughlin considered higher degree versions of the
tame symbol construction, \cite{brymcl:deg4_I,brymcl:deg4_II},
namely cup products of higher degree Deligne cohomology classes:
$\tame{f}{L}$ for $f$ a holomorphic invertible function and $L$ a
holomorphic line bundle, and $\tame{L}{L'}$ for a pair of
holomorphic line bundles. The geometric interpretation of the
symbols so obtained, also put forward in
refs.~\cite{brymcl:deg4_I,brymcl:deg4_II}, is that $\tame{f}{L}$
is a gerbe on $X$ with band ($\coin$ lien) $\sho{X}^\unit$ and a
holomorphic connective structure. A similar statement holds for
the $2$-gerbe $\tame{L}{L'}$.

\subsubsection{Cup products}
\label{sec:higher-cup-products}

From the point of view of cohomology classes, one computes the
relevant cup products.  Using~\eqref{eq:4}, we find that
\begin{math}
  \tame{f}{L}\in \delH[3]{X}{\ZZ}{2}
\end{math}
is represented by the cocycle
\begin{equation}
  \label{eq:39}
  \big( g_{jk}^{-m_{ij}}\,,
  -\frac{1}{\tate}\log_i\! f\: d\log g_{ij} \bigr)\,,
\end{equation}
having made the standard choices for $\log_i f$ and the
transition functions $g_{ij}$ of $L$ with respect to the choice
of a cover $\cover{U}_X$. Similarly, if $g'_{ij}$ are the
transition functions of $L'$, and $\tate c_{ijk}$ represents
$c_1(L)$ with respect to the cover $\cover{U}_X$, then
\begin{math}
  \tame{L}{L'}\in \delH[4]{X}{\ZZ}{2}
\end{math}
is represented by the cocycle
\begin{equation}
  \label{eq:40}
  \big( {g'_{kl}}^{-c_{ijk}}\,,
  -\frac{1}{\tate}\log g_{ij} \: d\log g'_{jk} \bigr)\,.
\end{equation}

\subsubsection{Hermitian variant}
\label{sec:high-herm-prod}

If we use the product
\begin{displaymath}
  \deligne{\ZZ}{1} \otimes \deligne{\ZZ}{1}
  \lto \dhh{1}
\end{displaymath}
introduced in sect.~\ref{sec:hermitian-structure}, for $f$, $L$
and $L'$ as above we have 
\begin{displaymath}
  \begin{aligned}
    \delH[1]{X}{\ZZ}{1}&\otimes \delH[2]{X}{\ZZ}{1}&
    \lto&    \dhhH[3]{X}{1}\\
    f&\otimes [L]& \longmapsto& \tamehh{f}{L}
  \end{aligned}
\end{displaymath}
Using the same \cech\ data as before, the symbol $\tamehh{f}{L}$
is represented by the cocycle
\begin{equation}
  \label{eq:41}
  \big( g_{jk}^{-m_{ij}}\,,
  -\frac{1}{\tate}\pi_1(\log_i\! f)\,\pi_0(\log g_{ij}) \bigr)\,.
\end{equation}
Similarly, with $L$ and $L'$ we have the product
\begin{displaymath}
  \begin{aligned}
    \delH[2]{X}{\ZZ}{1}&\otimes \delH[2]{X}{\ZZ}{1}&
    \lto&    \dhhH[4]{X}{1}\\
    [L]&\otimes [L']& \longmapsto& \tamehh{L}{L'}
  \end{aligned}
\end{displaymath}
and the representing cocycle
\begin{equation}
  \label{eq:42}
  \big( {g'_{kl}}^{-c_{ijk}}\,,
  -\frac{1}{\tate}\pi_1(\log g_{ij})\,\pi_0(\log g'_{jk})
  \bigr)\,.
\end{equation}
Similarly to the proof of prop.~\ref{prop:2}, the maps of
complexes 
\begin{math}
  \deligne{\ZZ}{2} \to \deligne{\ZZ}{1}
\end{math}
and 
\begin{math}
  \dhh{1} \to \deligne{\ZZ}{1}
\end{math}
induce corresponding maps on the symbols 
\begin{math}
  \tame{f}{L}
\end{math}
and 
\begin{math}
  \tamehh{f}{L}\,,
\end{math}
moreover their images agree in 
\begin{math}
  \delH[3]{X}{\ZZ}{1}\,.
\end{math}
An identical statement holds for
\begin{math}
  \tame{L}{L'}
\end{math}
and 
\begin{math}
  \tamehh{L}{L'}\,.
\end{math}

\subsection{Gerbes with Hermitian structure}
\label{sec:hermitian-structures}

Let $\gerbe{G}$ be a gerbe on $X$ with band $\sho{X}^\unit$
(\cite{MR49:8992}). After \cite{MR95m:18006,bry:loop}, its class
is an element of
\begin{math}
  \delH[3]{X}{\ZZ}{1}\iso H^2(X,\sho{X}^\unit)\,.
\end{math}
Let
\begin{math}
  \she[0]{X,+}
\end{math}
be the sheaf of real positive smooth functions on $X$.
\begin{definition}
  \label{def:2}
  A \emph{hermitian structure} on $\gerbe{G}$ consists of the
  following data:
  \begin{enumerate}
  \item\label{item:1} For each object $P$ in $\gerbe{G}_U$, is
    assigned a $\she[0]{U,+}$-torsor $\herm{P}$ (a
    $\RR_{+}$-principal bundle). The assignment must be
    compatible with the restriction functors
    \begin{math}
      i^*\colon \gerbe{G}_U\to    \gerbe{G}_V
    \end{math}
    arising from $i\colon V\hookrightarrow U$ in the cover
    $\cover{U}_X$ of $X$.
  \item\label{item:2} For each morphism $f\colon P\to Q$ in
    $\gerbe{G}_U$ a corresponding morphism
    \begin{math}
      f_*\colon \herm{P}\to \herm{Q}
    \end{math}
    of $\she[0]{U,+}$-torsors.\footnote{A $\she[0]{U,+}$-torsor
      will in general be automatically trivializable. However, in
      this context it is convenient to ``forget'' the actual
      trivializing map.} This map must be compatible with
    compositions of morphisms in $\gerbe{G}_U$ and with the
    restriction functors.
    
    For an object $P$ of $\gerbe{G}_U$, an automorphism
    $\varphi\in \aut{P}$ is identified with a section of
    $\sho{X}^\unit$ over $U$. We then require that
    \begin{equation}
      \label{eq:43}
      \begin{aligned}
        \varphi_*\colon  \herm{P} &\overset{\simeq}{\lto} \herm{P}\\
        h \;&\longmapsto h\cdot\abs{\varphi}^2
      \end{aligned}
    \end{equation}
    where the latter is the $\she[0]{U,+}$-action on the torsor
    $\herm{P}$.
  \end{enumerate}
\end{definition}
\begin{theorem}
  \label{thm:1}
  Equivalence classes of $\sho{X}^\unit$-gerbes with hermitian
  structure are classified by the group
  \begin{displaymath}
    \HHH^3 \bigl(X,
    \ZZ(1)_X \to \sho{X} \to \she[0]{X}\bigr)\,.
  \end{displaymath}
\end{theorem}
\begin{proof}
  Let $\gerbe{G}$ be an $\sho{X}^\unit$-gerbe on $X$ with
  hermitian structure as per definition~\ref{def:2}. Choose a
  full decomposition (see \cite{MR95m:18006}) with objects
  $P_i$ of $\gerbe{G}_{U_i}$ and isomorphisms
  \begin{math}
    f_{ij}\colon P_j\vert_{U_{ij}} \to P_i\vert_{U_{ij}}
  \end{math}
  with respect to a cover $\cover{U}_X$ of $X$. By a standard
  procedure (see refs.\cite{MR95m:18006,bry:loop}) these
  data determine a cochain
  \begin{math}
    g_{ijk}\in \aut{P_i}\vert_{U_{ijk}}
    \iso \sho{X}^\unit\vert_{U_{ijk}}
  \end{math}
  satisfying the cocycle condition and determining a class in
  $H^2(X,\sho{X}^\unit)$. Furthermore, choose sections $r_i$ of
  the torsors $\herm{P_i}$ above $U_i$. From
  condition~\ref{item:2} in definition~\ref{def:2} we have that
  there must exist $\rho_{ij} \in \she[0]{X,+}\vert_{U_{ij}}$
  such that:
  \begin{equation}
    \label{eq:44}
    {f_{ij}}_* (r_j) = r_i\cdot \rho_{ij}\,.
  \end{equation}
  On the $3$-skeleton of the cover we have that on one hand
  \begin{equation}
    \label{eq:45}
    {f_{ij}}_*\circ {f_{jk}}_* (r_k) = {f_{ij}}_*(r_j)\cdot \rho_{jk}
    =r_i\cdot \rho_{ij}\,\rho_{jk}\,,
  \end{equation}
  whereas on the other hand, since $f_{ij}\circ f_{jk} =
  g_{ijk}\circ f_{ik}\,,$ we have
  \begin{equation}
    \label{eq:46}
    (f_{ij}\circ f_{jk})_*(r_k)
    = {g_{ijk}}_*\circ {f_{ik}}_*(r_k)
    = {g_{ijk}}_*(r_i\cdot \rho_{ik})
    = r_i\cdot \abs{g_{ijk}}^2\rho_{ik}\,.
  \end{equation}
  Equating the right hand sides of eqs.~\eqref{eq:45}
  and~\eqref{eq:46}, and extracting the appropriate logarithms,
  we see we have obtained a \cech\ cocycle representing a class
  in 
  \begin{equation}
    \label{eq:47}
    \Check\HHH^3 \bigl(\cover{U}_X,
    \ZZ(1)_X \to \sho{X} \to \she[0]{X}\bigr)\,.
  \end{equation}
  
  Conversely, let a class in $\dhhH[3]{X}{1}$ be given, and
  assume we represent it via the choice of $\cover{U}_X$ by a
  degree $2$ \cech\ cocycle with values in the complex
  \begin{displaymath}
    \ZZ(1)_X \to \sho{X} \to \she[0]{X}\,,
  \end{displaymath}
  which we write as
  \begin{displaymath}
    \bigl(\tate c_{ijkl}\,,\:
    \log g_{ijk}\,,\:
    \onehalf\log \rho_{ij}\bigr)\,.
  \end{displaymath}
  This cocycle determines, via the map $\dhh{1} \to
  \deligne{\ZZ}{1}$, a cocycle $\{g_{ijk}\}\in
  \vC[2]{\cover{U}_X}{\sho{X}^\unit}$ which can be used,
  according to refs.~\cite{MR95m:18006,bry:loop}, to glue
  the local stacks $\tors{\sho{U_i}}$ into a global $\gerbe{G}$,
  in fact a gerbe. Given a $\sho{U_i}^\unit$-torsor $P_i$, namely
  an object of $\gerbe{G}_{U_i}\iso \tors{\sho{U_i}}$, define a
  hermitian structure by:
  \begin{displaymath}
    \herm{P_i} = \text{trivial}\;\she[0]{U_i,+}-\text{torsor}
  \end{displaymath}
  Then use $\rho_{ij}$ to glue $\herm{P_i}$ and $\herm{P_j}$ over
  $U_{ij}$, namely \emph{define} an isomorphism via
  eq.~\eqref{eq:44}. Since the isomorphisms $P_k\to P_i$ and
  $P_k\to P_j\to P_i$ differ by the equivalence determined by
  $g_{ijk}$, we see using~\eqref{eq:43} that the condition
  \begin{displaymath}
    \rho_{ij}\,\rho_{jk} = \abs{g_{ijk}}^2\rho_{ik}\,,
  \end{displaymath}
  ensuing from the cocycle condition, ensures the compatibility
  of this definition over $U_{ijk}$.
\end{proof}
\begin{corollary}
  Using the quasi-isomorphism
  \begin{displaymath}
    \dhh{1} \lqi
    \bigl(\ZZ(1)_X \to \sho{X} \to \she[0]{X}\bigr)\,,
  \end{displaymath}
  the class of a gerbe with hermitian structure is in fact in
  $\dhhH[3]{X}{1}\,.$
\end{corollary}
We will see (cf. sect.~\ref{sec:herm-conn-struct}) this group
also automatically classifies a special type of connective
structure on $\gerbe{G}$.

\subsection{Hermitian connective structure}
\label{sec:herm-conn-struct}

The structure defined in sect.~\ref{sec:hermitian-structures} can
be supplemented by a variant of Brylinski's connective
structure~\cite{bry:loop} by taking into account the first Hodge
filtration as in ref.~\cite{bry:quillen}. Let $\gerbe{G}$ be an
$\sho{X}^\unit$ gerbe over $X$.
\begin{definition}
  \label{def:3}
  A type $(1,0)$ \emph{connective structure} on $\gerbe{G}$ is
  the assignment to each object $P$ of $\gerbe{G}_U$ of a
  $F^1\!\sha[1]{U}$-torsor $\conn{P}$ compatible with restriction
  functors and morphisms of objects. In particular, for
  $\varphi\in \aut{P}$, we require that 
  \begin{equation}
    \label{eq:48}
    \begin{aligned}
      \varphi_*\colon  \conn{P} &\overset{\simeq}{\lto} \conn{P}\\
      \nabla \;&\longmapsto \nabla + d\log \varphi
    \end{aligned}
  \end{equation}
  where $\nabla$ is a section of $\conn{P}$ over
  $U$.\footnote{Note that $d\log\varphi$ is holomorphic, hence of
  type $(1,0)$.}
\end{definition}
\begin{definition}
  \label{def:4}
  Let $\gerbe{G}$ be equipped with a hermitian structure.  A type
  $(1,0)$ connective structure on $\gerbe{G}$ is
  \emph{compatible} with the hermitian structure if for each
  object $P$ of $\gerbe{G}$ there is an isomorphism of torsors
  \begin{align*}
    \herm{P} &\lto \conn{P}\\
    r &\longmapsto \nabla_r
  \end{align*}
  such that for a positive function $\rho$ on $U$
  \begin{displaymath}
    r\cdot\rho\longmapsto \nabla_r +\del\log\rho\,.
  \end{displaymath}
  (In other words, $\nabla_{r\cdot \rho}=\nabla_r
  +\del\log\rho\,.$)
\end{definition}
Connective structures of type $(1,0)$ are classified as follows.
\begin{theorem}
  \label{thm:2}
  Let again $\dhh{1}$ be the complex given by~\eqref{eq:15} for
  $l=1$. Equivalence classes of connective structures on a
  $\sho{X}^\unit$-gerbe $\gerbe{G}$ compatible with a given
  hermitian structure are classified by the group
  \begin{displaymath}
    \HHH^3 \bigl(X,\dhh{1}\bigr)\,.
  \end{displaymath}
\end{theorem}
We have the following analog of the existence and uniqueness of
the canonical connection on an invertible sheaf.
\begin{corollary}
  A connective structure compatible with a hermitian structure on
  a gerbe $\gerbe{G}$ is uniquely determined up to
  equivalence.
\end{corollary}
\begin{proof}
  It is an immediate consequence of the fact that the groups in
  Theorems~\ref{thm:1} and~\ref{thm:2}, being computed from
  quasi-isomorphic complexes, are actually the same (and equal to
  $\dhhH[3]{X}{1}$.)
\end{proof}
\begin{remark}
  \label{rem:3}
  The group $\HHH^3 \bigl(X,\dhh{1}\bigr)\iso \dhhH[3]{X}{1}$ is
  \emph{not} equal to Brylinski's
  \begin{equation*}
    \HHH^3\bigl(X,\ZZ(1)\to\she[0]{X}(1)\to\she[1]{X}(1)\bigr)\,,
  \end{equation*}
  cf.\ ref.\ \cite[Proposition 6.9 (1)]{bry:quillen}. (In fact
  there is an epimorphism
  \begin{math}
    C(1)^\bullet \to
  \bigl( \ZZ(1)\to\she[0]{X}(1)\to\she[1]{X}(1) \bigr)
  \end{math}
  with non-trivial kernel.)  It follows that the notion of
  ``hermitian gerbes with hermitian connective structure'' in
  loc.\ cit.\ is not identical to our notion of
  $\sho{X}^\unit$-gerbe with hermitian structure and compatible
  type $(1,0)$ connective structure.
\end{remark}
\begin{proof}[Proof of Theorem~\protect\ref{thm:2}]
  Choose a cover $\cover{U}_X$ as usual and let
  $(P_i,f_{ij},r_i)$ be a decomposition of $\gerbe{G}$ and its
  hermitian structure as in the proof of Theorem~\ref{thm:1}.
  
  If $\gerbe{G}$ has a compatible type $(1,0)$ connective
  structure, we have a map $\herm{G_{U_i}} \ni r_i\mapsto\nabla_i
  \in \herm{G_{U_i}}$. For every isomorphism $f_{ij}$ the
  compatibility condition from Definition~\ref{def:4} determines
  a form
  \begin{displaymath}
    \xi_{ij}=\del\log\rho_{ij} \in F^1\!\sha[1]{X}(U_{ij})
  \end{displaymath}
  satisfying the condition
  \begin{equation}
    \label{eq:49}
    \xi_{jk} -\xi_{ik} +\xi_{ij} = d\log g_{ijk}\,.
  \end{equation}
  The imaginary $2$-form $\eta_{ij} \eqdef \delb \xi_{ij} =
  \delb\del \log\rho_{ij}$ then is a cocycle with values in
  $F^1\!\sha[2]{X}\cap \she[2]{X}(1)$.
  
  Altogether, $g_{ijk}$, $\onehalf\log\rho_{ij}$, $\xi_{ij}$ and
  $\eta_{ij}$ determine a cocycle of total degree $3$ in the
  \cech\ resolution $\vC{\cover{U}_X}{\dhh{1}}$.
  
  Conversely, given a degree $3$ cocycle with values in
  $\dhh{1}$, a gerbe $\gerbe{G}$ with hermitian structure can be
  obtained by gluing trivial $\sho{U_i}^\unit$-torsors and
  $\she[0]{U_i,+}$ torsors as in Theorem~\ref{thm:1}.
  Furthermore, define a map by assigning the trivial
  $F^1\!\sha[1]{U_i}$-torsor to the trivial
  $\she[0]{U_i,+}$-torsor by 
  \begin{displaymath}
    r \longmapsto \nabla_r\coin \del\log r\,.
  \end{displaymath}
  Clearly, this defines a type $(1,0)$ connective structure
  compatible with the hermitian structure on $\gerbe{G}$.
\end{proof}
\begin{remark}
  \label{rem:1}
  Note the proof of Theorem~\ref{thm:2} that $d\eta_{ij} = 0$,
  hence we obtain a class
  \begin{displaymath}
    [\eta_{ij}] \in \HHH^3\bigl(X, F^1\!\sha{X}\cap
    \sigma^{2}\she{X}(1) \bigr)
  \end{displaymath}
  which can be associated to $\gerbe{G}$ via the obvious map
  \begin{displaymath}
    \dhh{1} \lto F^1\!\sha{X}\cap \sigma^{2}\she{X}(1)\,.
  \end{displaymath}
  This class plays the same role for $\gerbe{G}$ as the (global)
  imaginary form $c_1(\rho) = \delb\del \log \rho_i$ for a
  metrized line bundle $(L,\rho)$.
\end{remark}
\begin{remark}[\emph{Hermitian curving}]
  \label{rem:2}
  An equivalent degree $3$
  cocycle can be obtained by introducing the cochain $K_i\in
  \sha[1,1]{X}\cap \she[2]{X}(1) (U_i)$ of imaginary $2$-forms
  such that
  \begin{displaymath}
    \delb\del \log\rho_{ij} = K_j -K_i\,,
  \end{displaymath}
  and the imaginary $3$-form $\Omega_i\coin \Omega\vert_{U_i}$
  such that
  \begin{displaymath}
    dK_i = \Omega\vert_{U_i}\,,
  \end{displaymath}
  where $\Omega\in F^1\!A^3(X)\cap E^3(X)(1)$ (global
  sections). We can regard $K_i$ as the hermitian \emph{curving}
  and $\Omega$ as the hermitian \emph{$3$-curvature,}
  respectively, of the type $(1,0)$ hermitian connection.
\end{remark}

\subsection{The symbol $\tamehh{f}{L}$}
\label{sec:tame-symbol-tamehhfl}

Given an invertible function $f$ and a line bundle $L$ we have
seen there is a product $\tamehh{f}{L}\in \dhhH[3]{X}{1}$. We
briefly give a geometric construction of the corresponding
hermitian-holomorphic gerbe.

We need to recall from \cite{brymcl:deg4_II} the construction of
the gerbe $\gerbe{C}$ underlying $\tame{f}{L}$. $\gerbe{C}$ is
the stackification of the following pre-stack $\gerbe{C}^0$.  For
$U\hookrightarrow X$ objects of the category $\gerbe{C}^0_U$ are
non vanishing sections of $L\vert_U$. If $s\in L\vert_U$, and non
vanishing, it is denoted $\tame{f}{s}$ as an object of
$\gerbe{C}^0_U$. Given another non vanishing section $s'$ of $L$
over $U$, there is $g\in \sho{U}^\unit$ such that $s'= s g$.
Morphisms from $\tame{f}{s'}$ to $\tame{f}{s}$ are given by
sections of the Deligne torsor $\tame{f}{g}$ over $U$. For a
third non vanishing section $s''$, with $s'' = s' g' = s g g'$,
composition of morphisms in the category $\gerbe{C}^0_U$
corresponds to the $K$-theoretic property of the Deligne torsor:
\begin{displaymath}
  \tame{f}{gg'} \iso \tame{f}{g} \otimes \tame{f}{g'}\,.
\end{displaymath}
Given a trivialization of $L$ by a collection $\lbrace
s_i\rbrace$ relative to a cover $\cover{U}_X = \lbrace
U_i\rbrace_{i\in I}$, with transition functions $g_{ij}\in
\sho{X}^\unit (U_{ij})$, the objects 
\begin{math}
  \tame{f}{s_i}
\end{math}
and the morphisms
\begin{displaymath}
  \phi_{ij} = \lbrace \log_i f,g_{ij} \rbrace \colon
  \tame{f}{s_j} \to \tame{f}{s_i}
\end{displaymath}
provide a decomposition of $\gerbe{C}$ in the sense of
\cite{MR95m:18006}. It follows that the automorphisms
\begin{equation}
  \label{eq:50}
  h_{ijk} = \phi_{ij}\otimes \phi_{jk} \otimes \phi_{ik}^{-1} =
  g_{jk}^{-m_{ij}} \in \aut{\tame{f}{s_i}\vert_{U_{ijk}}}\iso
  \sho{X}^\unit (U_{ijk})
\end{equation}
represent the cohomology class of $\gerbe{C}$ in
$\delH[3]{X}{\ZZ}{1}\iso H^2(X,\sho{X}^\unit)$.

Now define a \emph{hermitian structure} on $\gerbe{C}$ as
follows. To an object $\tame{f}{s}$ of $\gerbe{C}_U$ we assign
\begin{equation}
  \label{eq:51}
  \tame{f}{s} \rightsquigarrow 
  \herm{\tame{f}{s}} = \text{trivial $\she[0]{U,+}$-torsor.}
\end{equation}
Then, given a morphism 
\begin{math}
  \tame{f}{g}\ni \phi \colon \tame{f}{s'} \to \tame{f}{s}
\end{math}
in $\gerbe{C}_U$, with $s' = sg$ as above, we use the hermitian
structure on the Deligne torsor underlying $\tame{f}{g}$ defined
in sect.~\ref{sec:hermitian-structure}, Proposition~\ref{prop:2}.
Namely
\begin{equation}
  \label{eq:52}
  \begin{aligned}
    \phi_*\colon \herm{\tame{f}{s'}} &\lto \herm{\tame{f}{s}} \\
    h \; &\longmapsto h \cdot \norm{\phi}^2
  \end{aligned}
\end{equation}
where $h$ is a local section of $\herm{\tame{f}{s'}}$, to be
identified with one of $\she[0]{U,+}$ and $\norm{\phi}$ is the
length of the non-vanishing section $\phi$. We have the following
analog of Proposition~\ref{prop:2}:
\begin{proposition}
  \label{prop:4}
  The class of the gerbe $\gerbe{C}$ underlying the symbol
  $\tame{f}{L}$ with hermitian structure defined by
  eqs.~\eqref{eq:51} and~\eqref{eq:52} is given by the product
  \begin{math}
    \tamehh{f}{L}
  \end{math}
  in the group
  \begin{math}
    \HHH^3\bigl(X,\ZZ(1)_X \to \sho{X} \to \she[0]{X}\bigr) \iso
    \dhhH[3]{X}{1}\,.
  \end{math}
\end{proposition}
\begin{proof}
  We need to find the class of the $\gerbe{C}$ as in the proof of
  Thm.~\ref{thm:1} and show it coincides with $\tamehh{f}{L}$ as
  computed in eq.~\eqref{eq:41}. To this end, let us use the
  decomposition of $\gerbe{C}$ given by the objects
  $\tame{f}{s_i}$ and morphisms
  \begin{math}
    \phi_{ij} = \lbrace \log_i f,g_{ij} \rbrace \colon
    \tame{f}{s_j} \to \tame{f}{s_i}
  \end{math}
  for non vanishing sections $s_i\in L\vert_{U_i}$, as before.
  The class of $\gerbe{C}$ (without extra structures) is
  represented by the cochain $g^{-m_{ij}}_{jk}$ already appearing
  in eq.~\eqref{eq:50}.
  
  Furthermore, in the hermitian Deligne torsor $\tame{f}{g_{ij}}$
  over $U_{ij}$ the logarithm of the length of the section
  \begin{math}
    \phi_{ij} = \lbrace \log_i f,g_{ij} \rbrace
  \end{math}
  is given by
  \begin{displaymath}
    \sigma_{ij} \coin
    \onehalf \log \norm{\phi_{ij}}^2 \coin
    \onehalf \log\rho_{ij} =
    -\frac{1}{\tate}\,\pi_1(\log_i\! f) \log\abs{g_{ij}}\,,
  \end{displaymath}
  cf.~eq.~\eqref{eq:27}. Thus we have found the total cocycle
  representing $\tamehh{f}{L}$ as in eq.~\eqref{eq:41}. Indeed,
  by computing the \cech\ coboundary we find
  \begin{displaymath}
    \sigma_{ij} -\sigma_{ik} +\sigma_{jk}
    = -m_{ij} \log \abs{g_{jk}}\,,
  \end{displaymath}
  as desired.
\end{proof}

\subsection{Hermitian 2-Gerbes}
\label{sec:herm-2-gerbes}

Let us briefly extend the considerations outlined in the previous
sections to $2$-gerbes over $X$ bound by $\sho{X}^\unit$.  (An
extended exposition of the local geometry of $2$-gerbes is to be
found in ref.~\cite{MR95m:18006}. See also \cite{brymcl:deg4_I}
for the abelian case.)

Recall that a $2$-gerbe $\twogerbe{G}$ over $X$ bound by a sheaf
of \emph{abelian groups} $\sheaf{H}$ is a fibered $2$-category
over $X$ which satisfies the $2$-descent condition for objects,
and such that for any two objects $P$ and $Q$ in the fiber
$2$-category $\twogerbe{G}_U$ over $U\subset X$ the fibered
category $\sheafhom (P,Q)$ is a stack. If fact, this fibered
category turns out to be an $\sheaf{H}$-gerbe equivalent to the
neutral one $\tors{\sheaf{H}}$. The properties of interest to us
are the following: $\twogerbe{G}$ is \emph{locally non-empty,}
namely there is a cover $\cover{U}_X$ of $X$ such that for
$U\subset X$ in the cover, the object set of $\twogerbe{G}_U$ is
non-empty; $\twogerbe{G}$ is \emph{locally connected,} namely any
two objects can be connected by a weakly invertible $1$-arrow
(that is, invertible up to a $2$-arrow); any two $1$-arrows can
be (locally) joined by a $2$-arrow; finally, for every $1$-arrow
its automorphism group is isomorphic in a specified way to
$\sheaf{H}$.

Once the appropriate notion of isomorphism for $2$-gerbes is
introduced, isomorphism classes of $2$-gerbes bound by
$\sheaf{H}$ are classified by the sheaf cohomology group
\begin{math}
  H^3(X,\sheaf{H})\,,
\end{math}
see, e.g. refs.~\cite{MR95m:18006,brymcl:deg4_I}.

In what follows, we shall set $\sheaf{H}=\sho{X}^\unit$. Hence we
can rephrase the previous statement by saying that isomorphism
classes of $2$-gerbes bound by $\sho{X}^\unit$ are classified by
the group
\begin{displaymath}
  H^3(X,\sho{X}^\unit) \iso \delH[4]{X}{\ZZ}{1}\,.
\end{displaymath}

We shall need the local calculation leading to the
classification, so we recall it here. Given a $2$-gerbe
$\twogerbe{G}$, let us choose a decomposition by selecting a
cover $\cover{U}_X$ of $X$ and a collection of objects $P_i$ in
$\twogerbe{G}_{U_i}$. There is a $1$-arrow
\begin{displaymath}
  f_{ij}\colon P_j\to P_i
\end{displaymath}
between their restrictions to $\twogerbe{G}_{U_{ij}}$.
Furthermore, from the axioms there is a $2$-arrow
\begin{displaymath}
  \alpha_{ijk} \colon f_{ij}\circ f_{jk} \Longrightarrow f_{ik}\,.
\end{displaymath}
Further restricting over a $4$-fold intersection
$U_{ijkl}$, we have two $1$-arrows
\begin{math}
  f_{ij}\circ f_{jk}\circ f_{kl}\colon P_l\to P_i
\end{math}
and
\begin{math}
  f_{il}\colon P_l\to P_i
\end{math}
and between them \emph{two} $2$-arrows, namely
\begin{math}
  \alpha_{ijl} \circ (\Id_{f_{ij}} * \alpha_{jkl})
\end{math}
and
\begin{math}
  \alpha_{ikl} \circ (\alpha_{ijk} * \Id_{f_{kl}})\,.
\end{math}
Since $2$-arrows are strictly invertible, it follows again from
the axioms that there exists a section $h_{ijkl}$ of
$\sho{X}^\unit$ over $U_{ijkl}$ such that
\begin{equation}
  \label{eq:57}
  \alpha_{ijl} \circ (\Id_{f_{ij}} * \alpha_{jkl})
  = h_{ijkl}\circ
  \alpha_{ikl} \circ (\alpha_{ijk} * \Id_{f_{kl}})\,.
\end{equation}
This section is a $3$-cocycle and the assignment $\twogerbe{G}
\mapsto [h]$ gives the classification isomorphism.

In analogy with what was previously done for gerbes, we are going
to define a notion of hermitian structure and of type $(1,0)$
\emph{connectivity} for $2$-gerbes on $X$ bound by
$\sho{X}^\unit$. Brylinski and McLaughlin defined a \emph{concept
  of connectivity} on a $2$-gerbe $\twogerbe{G}$ over $X$ to be
the datum of a compatible class of connective structures on the
gerbes
\begin{math}
  \sheafhom_U (P,Q)
\end{math}
for two objects $P$, $Q$ in the fiber $\twogerbe{G}_U$. It is
possible to introduce several variants of this notion, as done in
refs.~\cite{brymcl:deg4_I,brymcl:deg4_II}. Thus a type $(1,0)$
connectivity will just be the requirement that these connective
structures take their values in 
\begin{math}
  F^1\!\sha[1]{X}-\text{torsors.}
\end{math}

Let us model the concept of hermitian structure on a $2$-gerbe
after the one for gerbes given above in definition~\ref{def:2}. 
\begin{definition}
  \label{def:5}
  A \emph{hermitian structure} on a $\sho{X}^\unit$-$ 2$-gerbe
  $\twogerbe{G}$ over $X$ consists of the following data.
  \begin{enumerate}
  \item To each object $P$ in the fiber $2$-category
    $\twogerbe{G}_U$ over $U\subset X$ we assign a
    $\she[0]{U,+}$-gerbe $\herm{P}$ over $U$. (As before,
    $\she[0]{U,+}$ is the sheaf of real positive functions on
    $U$.) 
  \item This assignment must be compatible with the inverse image
    $2$-functors
    \begin{math}
      i^*\colon \twogerbe{G}_U \to \twogerbe{G}_V\,,
    \end{math}
    natural transformations 
    \begin{math}
      \varphi_{i,j}\colon j^*i^* \Rightarrow (ij)^*
    \end{math}
    and modifications 
    \begin{math}
      \alpha_{i,j,k} \colon
      \varphi_{ij,k} \circ (h^* * \varphi_{i,j})
      \Rrightarrow
      \varphi_{i,jk} \circ (\varphi_{j,k} * i^*)
    \end{math}
    arising from the inclusions
    \begin{math}
      i\colon V\hookrightarrow U,
    \end{math}
    \begin{math}
      j\colon W\hookrightarrow V,
    \end{math}
    and
    \begin{math}
      k\colon Z\hookrightarrow W,
    \end{math}
    in the cover $\cover{U}_X$.
  \item For each $1$-arrow $f\colon P\to Q$ in $\twogerbe{G}_U$ a
    corresponding equivalence 
    \begin{math}
      f_*\colon \herm{P}\to \herm{Q}
    \end{math}
    of $\she[0]{U,+}$-gerbes. For each $2$-arrow
    \begin{math}
      \alpha \colon f \Rightarrow f'
    \end{math}
    a corresponding natural transformation
    \begin{math}
      \alpha_*\colon f_*\Rightarrow f'_*
    \end{math}
    between equivalences. We ask that this correspondence be
    compatible with compositions of $1$- and $2$-arrows. Namely,
    for $1$-arrows $f,f'\colon P \to Q$ and $g,g'\colon Q\to R$
    and for $2$-arrows
    \begin{math}
      \alpha\colon f\Rightarrow f'
    \end{math}
    and
    \begin{math}
      \beta\colon g\Rightarrow g'
    \end{math}
    in $\twogerbe{G}_U$, which we compose as
    \begin{math}
      \beta *\alpha \colon g\circ f \Rightarrow g'\circ f',
    \end{math}
    we find a diagram of natural transformations
    \begin{equation}
      \label{eq:59}
      \xymatrix@!C{%
        g_*\circ f_* \ar@{=>}[d]_{\beta_* *\alpha_*}
        \ar@{=>}[r]^{\varepsilon(f,g)}
        & (g\circ f)_* \ar@{=>}[d]^{(\beta *\alpha)_*} \\
        g'_*\circ f'_* \ar@{=>}[r]_-{\varepsilon(f',g')}
        & (g'\circ f')_*
      }
    \end{equation}
    of equivalences between the $\she[0]{U,+}$-gerbes $\herm{P}$
    and $\herm{R}$ on $U\subset X.$
  \item From the axioms, the group of automorphisms of a
    $1$-arrow $f\colon P\to Q$ in $\twogerbe{G}_U$ is identified
    with $\sho{U}^\unit$. It follows that such an automorphism
    $\alpha$ (that is, a $2$-arrow from $f$ to itself) can be
    identified with a section $a\in \sho{U}^\unit$. We then
    require that the induced natural isomorphism
    \begin{displaymath}
      \alpha_* \colon f_* \Longrightarrow f_*\,, \quad
      \text{where} \;
      f_* \colon \herm{P} \lto \herm{Q}
    \end{displaymath}
    be identified with a section of $\she[0]{U,+}$ via the map
    \begin{equation}
      \label{eq:61}
      a \longmapsto \abs{a}^2
    \end{equation}
    and an appropriate labeling of $\herm{P}$ and $\herm{Q}$ by
    objects $r$ and $s$, respectively. In more detail, given an
    arrow $f_*(r) \to s$ in $\herm{Q}$, the action of $\alpha$
    via $\alpha_*$ will amount to an automorphism of $s$. We
    require that it be $\abs{a}^2$.
  \end{enumerate}
\end{definition}
\begin{remark}
  The abstract nonsense of definition~\ref{def:5} could have more
  succinctly characterized by saying that the correspondence
  $\herm{\cdot}$ realizes a Cartesian $2$-functor between
  $\twogerbe{G}$ and the $2$-gerbe
  $\gerbes{\she[0]{X,+}}$ on $X$, shifting to the
  reader the burden of unraveling the diagrams.
\end{remark}
We have the following analog of theorem~\ref{thm:1}:
\begin{theorem}
  \label{thm:3}
  Isomorphism classes of $\sho{X}^\unit$-$2$-gerbes with
  hermitian structure in the sense of definition~\ref{def:5} are
  classified by the group
  \begin{displaymath}
    \HHH^4 \bigl(X,
    \ZZ(1)_X \to \sho{X} \to \she[0]{X}\bigr)\iso \dhhH[4]{X}{1}\,.
  \end{displaymath}
\end{theorem}
\begin{proof}
  Let $\twogerbe{G}$ be a $\sho{X}^\unit$-$2$-gerbe on $X$ with
  hermitian structure as per definition~\ref{def:5}. Forgetting
  the hermitian structure, $\twogerbe{G}$ will determine a class
  in the group
  \begin{math}
    \delH[4]{X}{\ZZ}{1} \iso H^3(X,\sho{X}^\unit),
  \end{math}
  and we have briefly recalled before --- cf. eq.~\eqref{eq:57}
  --- how to obtain a $3$-cocycle representing the class
  of~$\twogerbe{G}$.
  
  To obtain the rest of the cocycle with values in the complex
  \begin{math}
    \ZZ(1)_X \to \sho{X} \to \she[0]{X}
  \end{math}
  let us make the same choice for a decomposition of
  $\twogerbe{G}$ with respect to the cover $\cover{U}_X$: a
  collection of objects $P_i$ in $\twogerbe{G}_{U_i}$, $1$-arrows
  $f_{ij} \colon P_j\to P_i$ between their restrictions and
  $2$-arrows
  \begin{math}
    \alpha_{ijk}\colon f_{ij}\circ f_{jk} \Rightarrow f_{ik}\,.
  \end{math}
    
  We shall also need a decomposition of the
  $\she[0]{U_i,+}$-gerbes $\herm{P_i}$: to this end let us choose
  objects $r_i$ over $U_i$ and arrows
  \begin{math}
    \xi_{ij}\colon (f_{ij})_*(r_j)\to r_i
  \end{math}
  between their restriction to $U_{ij}$.
  
  Let us consider a triple of objects $P_i, P_j, P_k$ over
  $U_{ijk}$. (we are implicitly restricting to the fiber
  $2$-category $\twogerbe{G}_{U_{ijk}}$.) We obtain the following
  diagram in $\herm{P_i}\rvert_{U_{ijk}}$:
  \begin{equation}
    \label{eq:58}
    \xymatrix@!C{%
      (f_{ij})_* (f_{jk})_* (r_k) \ar[d]
      \ar[r]^-{(f_{ij})_* (\xi_{jk})}
      & (f_{ij})_*(r_j) \ar[d]^{\xi_{ij}} \\
      (f_{ik})_* (r_k) \ar[r]_{\xi_{ik}}
      & r_i \ar @(dr,ur)_{\rho_{ijk}}
    }
  \end{equation}
  The left vertical arrow in~\eqref{eq:58} results from the
  composition of two-arrows
  \begin{displaymath}
    \xymatrix@!C{%
      (f_{ij})_* \circ (f_{jk})_*
      \ar@{=>}[r]^{\varepsilon_{ijk}}
      & (f_{ij}\circ f_{jk})_*
      \ar@{=>}[r]^{(\alpha_{ijk})_*}
      & (f_{ik})_*}
  \end{displaymath}
  resulting from diagram~\eqref{eq:59} in definition~\ref{def:5}.
  At the level of objects in the gerbe $\herm{P_i}$
  diagram~\eqref{eq:59} is of course not commutative, so we
  obtain a section $\rho_{ijk}\in \aut{r_i}$, which we can
  identify with a section of the sheaf $\she[0]{U,+}$ over
  $U_{ijk}$.
  
  Now consider a four-fold intersection $U_{ijkl}$: we have a
  cube determined by the objects $r_i,\dots, r_l$ whose faces are
  built from copies of \eqref{eq:58}. Since this cube brings in
  the relation~\eqref{eq:57}, using the mapping of the
  $\sho{X}^\unit$ action spelled out in the last point in
  definition~\ref{def:5}, we get the relation
  \begin{equation}
    \label{eq:60}
    \rho^{\phantom{-1}}_{jkl}\,
    \rho^{-1}_{ikl}\,
    \rho^{\phantom{-1}}_{ijl}\,
    \rho^{-1}_{ijk}
    = \abs{h_{ijkl}}^2
  \end{equation}
  which, after taking the appropriate logarithms, defines a
  \cech\ cocycle representing a class in
  \begin{displaymath}
    \Check\HHH ^4\bigl(\cover{U}_X, \ZZ(1)_X \to \sho{X} \to
    \she[0]{X}\bigr)\,.
  \end{displaymath}
  Details (and diagram chasing) are straightforward and left to
  the reader.

  Conversely, let us be given a class in 
  \begin{displaymath}
    \HHH^4 \bigl(X,
    \ZZ(1)_X \to \sho{X} \to \she[0]{X}\bigr)
    \iso
    \HHH^3 \bigl(X,
    \sho{X}^\unit \overset{\abs{\cdot}}{\to} \she[0]{X,+}
    \bigr) \,,
  \end{displaymath}
  and let us assume it is represented by the (multiplicative)
  \cech\ cocycle
  \begin{math}
    \bigl(h_{ijkl},\rho_{ijk}\bigr)\,.
  \end{math}
  Let just explain the construction of a corresponding $2$-gerbe
  with hermitian structure (up to equivalence). Again, details
  will be left to the reader.

  We first apply the map
  \begin{displaymath}
    \bigl(\ZZ(1)_X \to \sho{X} \to \she[0]{X}\bigr)
    \lto \bigl( \ZZ(1)_X \to \sho{X} \bigr)
  \end{displaymath}
  to the representative \cech\ cocycle to reconstruct a
  $\sho{X}^\unit$-$2$-gerbe $\twogerbe{G}$ according to
  refs.~\cite{MR95m:18006,brymcl:deg4_I,brymcl:deg4_II}. Recall
  that this is accomplished by gluing the local stacks
  \begin{math}
    \gerbes{\sho{U_i}^\unit}
  \end{math}
  using $h_{ijkl}$. Secondly, we define a hermitian structure as
  follows.  Assign to any object $P_i$ over $U_i$ of the
  so-determined $2$-gerbe $\twogerbe{G}$ the trivial
  $\she[0]{U_i,+}$-gerbe $\herm{P_i}=\tors{\she[0]{U_i,+}}$. For
  a triple of such on $U_{ijk}$ we use $\rho_{ijk}\in
  \she[0]{U_i,+}\vert_{U_{ijk}}$ as an automorphism of an object
  $r_i$ in $\herm{P_i}$.

  Checking that this structure satisfies the properties in
  definition~\ref{def:5} and it defines a $2$-gerbe with
  hermitian structure whose class is the one we started with is
  modeled after the pattern of refs.~\cite{MR95m:18006}
  and~\cite{bry:loop} and it will be left to the reader.
\end{proof}

As mentioned before, a connectivity on a
$\sho{X}^\unit$-$2$-gerbe is in practice the assignment of
compatible connective structures on the local gerbes of
morphisms. We have the following definition (see also
\cite[sect.~7]{bry:quillen}, for the first part):
\begin{definition}
  \label{def:6}
  Let $\twogerbe{G}$ be a $\sho{X}^\unit$-$2$-gerbe on $X$.
  \begin{enumerate}
  \item A type $(1,0)$ \emph{concept of connectivity} on
    $\twogerbe{G}$ is the assignment of a $F^1\!\sha[1]{U}$-gerbe
    $\conn{P}$ to each object $P$ in $\twogerbe{G}_U$. This
    assignment will have to satisfy properties analogous to those
    of definition~\ref{def:5}. Of course, in the last condition,
    the map~\eqref{eq:61} will have to be replaced by
    \begin{math}
      a\mapsto d\log a\,.
    \end{math}
  \item A type $(1,0)$ concept of connectivity is
    \emph{compatible} with a hermitian structure if for each
    object $P$ of $\twogerbe{G}_U$ there is an equivalence of
    gerbes
    \begin{displaymath}
      \herm{P} \lto \conn{P}
    \end{displaymath}
    satisfying the obvious compatibility conditions with the
    operations of $\twogerbe{G}_U$ and the restrictions.
  \end{enumerate}
\end{definition}
The proof of the following theorem can be patterned after an
appropriate generalization of the proof of Theorem~\ref{thm:2},
so we shall omit it.
\begin{theorem}
  \label{thm:4}
  Let $\twogerbe{G}$ be a $\sho{X}^\unit$-$2$-gerbe with
  hermitian structure and let $\dhh{1}$ be the complex given
  by~\eqref{eq:15} for $l=1$. Equivalence classes of type
  $(1,0)$ connectivities on $\twogerbe{G}$ compatible with the
  given hermitian structure are classified by the group
  \begin{displaymath}
    \HHH^4 \bigl(X,\dhh{1}\bigr)\,.
  \end{displaymath}
  Furthermore, the equivalence class is unique.
\end{theorem}

\subsection{The symbol $\tamehh{L}{L'\,}$}
\label{sec:symbol-tamehhll}

We have seen that given two line bundles $L$ and $L'$ over $X$
their cup product $\tamehh{L}{L'}$ defines a class in
$\dhhH[4]{X}{1}$. According to Theorem~\ref{thm:3} it corresponds
to an equivalence class of $2$-gerbes with hermitian structure.
Using the obvious maps of complexes $\dhh{1}\to \deligne{\ZZ}{1}$
and $\deligne{\ZZ}{2}\to \deligne{\ZZ}{1}$, the geometric
$2$-gerbe $\twogerbe{G}$ that underlies $\tamehh{L}{L'}$ is the
same one as for the standard symbol $\tame{L}{L'}$ constructed by
Brylinski and McLaughlin.

Recall (see ref.~\cite{brymcl:deg4_II} for more details) that
objects of $\twogerbe{G}$ underlying $\tame{L}{L'}$ over
$U\subset X$ are the non-vanishing sections $s$ of $L\vert_U$,
denoted $\tame{s}{L}$. Given another non vanishing section $s'\in
L\lvert_U$ we have $s'=sg$ for an invertible function $g$ over
$U$. Then the category of morphisms from $\tame{s'}{L}$ to
$\tame{s}{L}$ is the \emph{gerbe} $\tame{g}{L}$ defined in
section~\ref{sec:tame-symbol-tamehhfl}. For a third non vanishing
section $s''$ of $L$ over $U$, with $s''=s'\, g'$, the morphism
composition functor is given by the equivalence
\begin{displaymath}
  \tame{g}{L'}\otimes \tame{g'}{L} \lto \tame{gg'}{L}
\end{displaymath}
where on the left hand side we have the contracted product of two
(abelian) gerbes. To be precise, it turns out that $\twogerbe{G}$
is an appropriate ``$2$-stackification'' of the $2$-pre-stack
defined here.

A calculation in ref.~\cite{brymcl:deg4_II} shows that with
respect to the trivializations $\lbrace g_{ij}\rbrace$ and
$\lbrace g'_{ij}\rbrace$ of $L$ and $L'$, respectively, the class
of $\twogerbe{G}$ is represented by the cocycle
\begin{math}
  {g'_{kl}}^{-c_{ijk}}\in \sho{X}^\unit(U_{ijkl})\,,
\end{math}
where the cocycle $c_{ijk}$ represents $c_1(L)$.

We can define a hermitian structure on $\twogerbe{G}$ as
follows.  To an object $\tame{s}{L'}$ of $\twogerbe{G}_U$ we
assign 
\begin{equation}
  \label{eq:62}
  \tame{s}{L'} \rightsquigarrow
  \herm{\tame{s}{L'}} = \text{trivial $\she[0]{U,+}$-gerbe.}
\end{equation}
Furthermore, as remarked above we have
\begin{math}
  \sheafhom_U \bigl(\tame{s'}{L'},\tame{s}{L'}\bigr)
  \iso \tame{g}{L'}\,.
\end{math}
Thus we set
\begin{equation}
  \label{eq:63}
  \sheafhom_U \bigl(
  \herm{\tame{s'}{L'}},\herm{\tame{s}{L'}}
  \bigr)
  = \tamehh{g}{L'}\,,
\end{equation}
where on the right hand side we use the hermitian structure on
the gerbe $\tame{g}{L'}$ as defined in
section~\ref{sec:tame-symbol-tamehhfl}. On the left hand side
of~\eqref{eq:63} we have the equivalences of the two
$\she[0]{U,+}$-gerbes.

The proof of the following proposition is a straightforward
generalization of the one for proposition~\ref{prop:4}.
\begin{proposition}
  The class of the $\sho{X}^\unit$-$2$-gerbe $\twogerbe{G}$
  underlying the symbol $\tame{L}{L'}$ with hermitian structure
  defined by eqs.~\eqref{eq:62} and~\eqref{eq:63} is given by
  the product $\tamehh{L}{L'}$ in the group
  \begin{math}
    \HHH^4\bigl(X,\ZZ(1)_X \to \sho{X} \to \she[0]{X}\bigr) \iso
    \dhhH[4]{X}{1}\,.
  \end{math}
\end{proposition}

\subsection{Comparisons and relations with other definitions}
\label{sec:comp-relat-with}

Recall from refs.\ \cite{brymcl:deg4_I,brymcl:deg4_II}, that
analytic connective structures on gerbes with band
$\sho{X}^\unit$ are classified by the group
$\delH[3]{X}{\ZZ}{2}$. Similarly, for $2$-gerbes with the same
band, the relevant group is $\delH[4]{X}{\ZZ}{2}$. In the
previous sections we have introduced hermitian structures and
type-$(1,0)$ connective structures on gerbes and $2$-gerbes with
band $\sho{X}^\unit$. We define the concept of compatibility
analogously to the case of line bundles in sect.\
\ref{sec:comparisons} as follows.

Let $\gerbe{G}$ be a $\sho{X}^\unit$-gerbe on $X$. Let
$\conn{\cdot}^\mathit{an}$ be a (holomorphic) connective structure on
$\gerbe{G}$ in the sense of refs.\
\cite{brymcl:deg4_I,brymcl:deg4_II}, and let $\conn{\cdot}^h$ be a
connective structure \emph{on the same gerbe} in the sense of
sect.\ \ref{sec:herm-conn-struct}.

The relevant group classifying $\gerbe{G}$ equipped with both
types of connections is therefore $\HHH^3(X,
\Tilde{\Gamma}(2)^\bullet)$, where the complex
$\Tilde{\Gamma}(2)^\bullet$ has been introduced in sect.\
\ref{sec:comparisons}.
\begin{definition}
  \label{def:7}
  We say that $\conn{\cdot}^\mathit{an}$ and
  $\conn{\cdot}^\mathit{h}$ are \emph{compatible} if for any
  object $P$ of $\gerbe{G}_U$, $U\subset X$, there is an
  isomorphism of torsors $\conn{P}^\mathit{an}\iso
  \conn{P}^\mathit{h}$ (after lambda-extension of
  $\conn{P}^\mathit{an}$ from $\shomega[1]{U}$ to
  $\sha[1,0]{U}$.)
\end{definition}
Similarly, if $\twogerbe{G}$ is a $\sho{X}^\unit$-$2$-gerbe on
$X$, carrying both types of connective structures, its class is
an element of the group $\HHH^4(X,
\Tilde{\Gamma}(2)^\bullet)$. We can also repeat the above
definition, taking care that now for any object of $\twogerbe{G}$
over $U\subset X$, $\conn{P}^\mathit{an}\iso \conn{P}^\mathit{h}$
must be an equivalence of gerbes.

The next lemma immediately follows from the definitions.
\begin{lemma} Let $\Gamma(2)^\bullet$ be the complex defined in
  sect.\ \ref{sec:comparisons}.
  \begin{enumerate}
  \item Classes of $\sho{X}^\unit$-gerbes with compatible
    connective structures in the sense of definition \ref{def:7}
    are classified by the elements of the group $\HHH^3(X, \Gamma
    (2)^\bullet)$.
  \item Similarly, classes of $\sho{X}^\unit$-$2$-gerbes with
    compatible connective structures are classified by $\HHH^4(X,
    \Gamma (2)^\bullet)$.
  \end{enumerate}
\end{lemma}

\subsubsection{Compatibility and flatness conditions}
\label{sec:flatness-conditions}

While these definitions seem to follow the pattern of line
bundles analyzed in sect.\ \ref{sec:comparisons}, there in an
important difference, namely gerbes (or $2$-gerbes) satisfying
the compatibility condition of definition~\ref{def:7} are
\emph{not necessarily flat!} Moreover, in the present framework
the compatibility condition is less special than it was seen in
the case of line bundles.  This is can be seen by way of the
following cohomological argument.

The complex $\Gamma (2)^\bullet$ introduced in sect.\
\ref{sec:comparisons} is easily seen to be a quotient of the
complex $\dhh{2}$:
\begin{equation*}
  \dhh{2} \lto \Gamma (2)^\bullet \lto 0\,.
\end{equation*}
The kernel is complicated, but up to quasi-isomorphism, it can be
reduced (by direct computation) to the one-element complex
$\she[2]{X}(1)\cap \sha[1,1]{X}[-4]$ so that we have the
triangle:
\begin{equation*}
  \she[2]{X}(1)\cap \sha[1,1]{X}[-4] \lto
  \dhh{2} \lto \Gamma (2)^\bullet \overset{+1}{\lto}
\end{equation*}
Focusing our attention to degree $3$ and $4$, we get the  sequence:
\begin{equation*}
  0 \to H^2 (X, \RR(2)/\ZZ(2))\to
  \HHH^3 (X, \Gamma (2)^\bullet)
  \to E^2(X)(1)\cap A^{1,1}(X) \to \dhhH[4]{X}{2} \to
  \HHH^4 (X, \Gamma (2)^\bullet)\to 0\,,
\end{equation*}
where we have used lemma~\ref{lem:2}.  Moreover, the exact
sequence from the proof of lemma~\ref{lem:4} relating
$\Tilde{\Gamma}(2)^\bullet$ to $\Gamma (2)^\bullet$ yields the
following completion of~\eqref{eq:64}:
\begin{gather*}
  0\to H^1(X,\RR(2)/\ZZ(2)) \to 
  \pic (X,\nabla,h) \to E^1(X)(1) \to
  \HHH^3 (X, \Gamma (2)^\bullet) \to
  \HHH^3 (X, \Tilde{\Gamma} (2)^\bullet) \to 0\\
  \intertext{and}
  \HHH^4 (X, \Gamma (2)^\bullet) \xrightarrow{\iso}
  \HHH^4 (X, \Tilde{\Gamma} (2)^\bullet)\,,
\end{gather*}
where we have used that $\she[1]{X}(1)$ is soft. In summary we have:
\begin{proposition}
  \label{prop:6}
  \hfill
  \begin{enumerate}
  \item The class of a $\sho{X}^\unit$-gerbe supporting both
    types of connective structures can be lifted to a class of
    compatible connective structures on a (possibly equivalent)
    gerbe.
  \item A $\sho{X}^\unit$-gerbe with compatible connective
    structures is flat if the (trivial) $(1,1)$-curving is zero
    (cf.\ sect.\ \ref{sec:herm-conn-struct}, remarks~\ref{rem:1}
    and~\ref{rem:2}.)
  \item A $\sho{X}^\unit$-$2$-gerbe supporting both types of
    connective structures is equivalent to a $2$-gerbe with
    compatible connective structures. Its class can be lifted to
    $\dhhH[4]{X}{2}$
  \end{enumerate}
\end{proposition}

\subsubsection{Comparing $\protect\tame{f}{L}$ and
  $\protect\tame{L}{L'}$ with their hermitian variants}
\label{sec:comp-prot-prot}

The higher symbols $\tame{f}{L}$ and $\tamehh{f}{L}$ have the
same underlying gerbe, and similarly $\tame{L}{L'}$ and
$\tamehh{L}{L'}$ determine the same $2$-gerbe. Let us denote
them, respectively, by $\tameg{f}{L}$ and
$\tameg{L}{L'}$. By construction, they determine classes in
$\HHH^3(X,\Tilde{\Gamma}(2)^\bullet)$ and
$\HHH^4(X,\Tilde{\Gamma}(2)^\bullet)$, respectively.
The proposition specializes to this case as follows:
\begin{corollary}
  The connective structures $\conn{\cdot}^{\mathit{an}}$ and
  $\conn{\cdot}^{\mathit{h}}$ on $\tameg{f}{L}$ are compatible
  (up to $\she[1]{U}$-torsor automorphism).

  The analytic and hermitian connective structures on the
  $2$-gerbe $\tameg{L}{L'}$ are compatible.
\end{corollary}
\begin{proof}
  The statement follows at once from the calculations preceding
  the proposition.
\end{proof}
\begin{remark}
As an alternative proof of the corollary, note that a calculation
analogous to that of the proof of proposition~\ref{prop:5}
from the cocycle representations~\eqref{eq:39} and~\eqref{eq:41},
yields the $1$-cocycle
\begin{math}
  r_2(f,g_{ij})
\end{math}
with values in $\she[1]{X}(1)$, where $g_{ij}$ are the transition
functions of $L$. This cocycle represents the zero class
(softness of $\she[1]{X}(1)$), therefore
\begin{math}
  r_2(f,g_{ij}) = \eta_j-\eta_i\,,
\end{math}
and this choice is determined up to a global section of
$\she[1]{X}(1)$.

Similarly, in the case of $\tameg{L}{L'}$ we get the $2$-cocycle
\begin{math}
  r_2(g_{ij},g'_{jk})
\end{math}
which again represents the zero class.
\end{remark}

\section{Concluding remarks}
\label{sec:conclusions}

In this paper we have put forward a definition for the concept of
hermitian structure, and associated compatible connective
structure for gerbes and $2$-gerbes with band $\sho{X}^\unit$. We
have presented classification results in terms of low degree
hermitian holomorphic Deligne cohomology groups. Notable examples
are provided by higher versions of the classical notion of tame
symbol associated to two invertible functions. Indeed, our second
main result that there exists a modified version of the cup
product in low degree Deligne cohomology taking values in the
first hermitian holomorphic Deligne complex, naturally provides
the symbols $\tame{f}{L}$ and $\tame{L}{L'}$ with hermitian
structures according to our definition.

Two questions naturally arise. Since $\tame{f}{L}$ and
$\tame{L}{L'}$ also carry an analytic connective structure, we
may ask to what degree the latter and the hermitian one are
compatible. Remark~\ref{rem:3} prompts a second obvious question
regarding the relation between our classification
theorems~\ref{thm:2} and~\ref{thm:4} and others', notably
Brylinski's (\cite[Proposition 6.9 (1)]{bry:quillen}).

We have analyzed the compatibility in cohomological terms, first
for line bundles (in the sense of $\sho{X}^\unit$-torsors) and
then for gerbes and $2$-gerbes with band-$\sho{X}^\unit$, with
somewhat surprising results. Whereas the compatibility may be
regarded as exceptional for a line bundle---and it implies its
flatness---it is not so for gerbes (or $2$-gerbes). Thus flatness
is not a necessary condition. In the specific case of the tame
symbols and their generalizations, we have found that while the
compatibility of $\tame{f}{g}$ and $\tamehh{f}{g}$ (that is,
their respective connections) may in general be obstructed,
$\tame{f}{L}$ and $\tamehh{f}{L}$ can always be made compatible,
and $\tame{L}{L'}$ and $\tamehh{L}{L'}$ are automatically so.

As for the relation with other notions of ``hermitian gerbe''
with ``hermitian connective structure'' (or $2$-gerbe) there
appear to be subtle differences in the definitions which we can
trace to what aspect of line bundles with connection we decide to
generalize. Our approach has been to copy the concept of
\emph{metrized analytic (or algebraic) line bundle} familiar from
Arakelov geometry (cf.\ ref.\ \cite{lang:arakelov}). On the other
hand, one could describe a metrized $\sho{X}^\unit$-line bundle
by means of the $\TT$-reduction of its associated smooth line
bundle plus a unitary connection. Whereas these two approaches
are equivalent in the case of line bundles, they seem to diverge
as soon as we move on to gerbes. (And possibly matters worsen in
the case of $2$-gerbes.) This may also serve to explain the lack
of uniqueness found by Hitchin's student D.~Chatterjee in his
thesis. Although that school's approach to gerbes lacks the
categorical input (in fact for them a gerbe is just the ``torsor
cocycle'' in the sense of \cite{MR95m:18006}) the definition of
hermitian gerbe is along Brylinski's lines.

Another difference is the following.  Our cohomological
characterization via the group $\dhhH[k]{X}{1}\iso \HHH^k(X,
\dhh{1})$, $k=3,4$, involves forms of degree two, which points to
a natural notion of curving naturally associated with the
structures we have defined (cf.\ remarks~\ref{rem:1}
and~\ref{rem:2}). This is obviously absent in the truncated group
in remark~\ref{rem:3}. The cohomological analysis of
sect.~\ref{sec:comp-relat-with}, where the group $\dhhH[4]{X}{2}$
appears, suggests that curvings can be a very nuanced structure,
however dealing with them in detail falls outside the scope of
the present work.

We hope to further elucidate matters in the future in another
publication.

\appendix

\section{Remarks on Hodge-Tate structures}
\label{sec:remarks-hodge-tate}

The relation between the ``imaginary part'' map made in
sect.~\ref{sec:hermitian-structure} together with the product
\begin{math}
  \deligne{\ZZ}{1}\otimes \deligne{\ZZ}{1}
  \to \tate\otimes\dhh{1}
\end{math},
and the cup product
\begin{math}
  \deligne{\ZZ}{1}\otimes \deligne{\ZZ}{1}
  \to \deligne{\ZZ}{2}
\end{math}
giving rise to the tame symbol becomes more transparent from the
point of view of Hodge-Tate structures.

\subsection{A Mixed Hodge Structure}
\label{sec:mixed-hodge-struct}

Let us briefly recall the following well known MHS on $\CC^3$,
see \cite{del:symbole,MR95a:19008}. Consider, as before,
\begin{equation}
  \label{eq:29}
  M^{(2)} = 
  \begin{pmatrix}
    1&& \\ x&1&\\ z& y & 1
  \end{pmatrix}
\end{equation}
with complex entries $x,y,z$. Consider also its canonical version
\begin{equation}
  \label{eq:30}
  A^{(2)} = 
  \begin{pmatrix}
    1&& \\ x&\tate &\\ z& \tate\, y & (\tate)^2
  \end{pmatrix}\,.
\end{equation}
The MHS $\hM_{2}$ corresponding to $M^{(2)}$, or more precisely
$A^{(2)}$, comprises the following data. The integer lattice is
the $\ZZ$ span of the columns of $A^{(2)}$, and similarly for
$\QQ$ and $\RR$.  Let $v_0, v_1, v_2$ denote the columns of
$A^{(2)}$ starting from the left. The weight spaces are
\begin{math}
  W_{-2k}\hM^{(2)} = \mathrm{span} \langle v_k,\dotsc,v_2\rangle
\end{math}
(over the appropriate ring), and the Hodge filtration is given by
\begin{math}
  F^{-k}\hM^{(2)}(\CC) = \CC
  \langle e_0,\dotsc,e_k\rangle\,,
\end{math}
where the $e_i$'s are the standard basis vectors in $\CC^2$. The
graded quotients
\begin{math}
  \gr^{W}_{-2k}\hM^{(2)}
\end{math}
are the Tate structures $\ZZ(0)$, $\ZZ(1)$, and $\ZZ(2)$. A
change of the generators $v_i$ preserving the structure clearly
amounts to a change of $A^{(2)}$ by right multiplication by a
lower unipotent matrix over $\ZZ$ (or $\QQ$ or $\RR$). This is
the same as changing $M^{(2)}$ by a matrix in $H_\ZZ$ (or the
appropriate ring thereof) as in sect.~\ref{sec:heisenberg-group}.
\footnote{These data correspond to the case $N=2$ of a MHS on
  $\CC^N$ defined for any integer $N$, cf.  \cite{MR95a:19008}}

The real structure underlying $\hM^{(2)}$ is linked to the
hermitian structure on the bundle $H_\CC/H_\ZZ$ as presented in
sect.~\ref{sec:remarks-heis-group}. In \cite{MR95a:19008} the
image of $A^{(2)}$ in $\GL_2(\CC)/\GL_2(\RR)$ is obtained by
computing the matrix
\begin{displaymath}
  B\eqdef A \Bar{A}^{-1}
  \left(
  \begin{smallmatrix}
    1 && \\ & -1 & \\ && 1
  \end{smallmatrix}
  \right)\,,
\end{displaymath}
(we have dropped the superscript $(2)$ for ease of notation).
The logarithm is:
\begin{displaymath}
  \onehalf\log B = 
  \begin{pmatrix}
    1 &&\\ \pi_0(x) & 1 & \\
    \pi_1(z) -\pi_1(x)\pi_0(y)&\pi_0(y) &1
  \end{pmatrix}\,.
\end{displaymath}
We immediately recognize the expression of the hermitian form as
given in sect.~\ref{sec:remarks-heis-group}.

\subsection{The big period }
\label{sec:big-period}
In ref. \cite{MR99i:19004} Goncharov defines a tensor
\begin{equation*}
  P(\hM) \in \CC\otimes_\QQ\CC
\end{equation*}
associated to a MHS (technically, a framed one) $\hM$.  For the
MHS defined by the period matrix~\eqref{eq:29} it is computed as
follows. Let $f_0, f_1, f_2$ be the dual basis to $v_0,v_1, v_2$.
Then, according to ref. \cite{MR99i:19004},
\begin{displaymath}
  P(\hM^{(2)}) = \sum_{k} \langle f_2, M^{(2)} v_k \rangle
  \otimes_\QQ
  \langle f_k, {M^{(2)}}^{-1}v_0\rangle \,.
\end{displaymath}
Performing the calculation we find:
\begin{equation}
  \label{eq:31}
  \begin{split}
    P(\hM^{(2)}) &= \frac{z}{(\tate)^2}\otimes 1
    -1\otimes \frac{z}{(\tate)^2}\\
    &+1\otimes \frac{xy}{(\tate)^2}
    -\frac{y}{\tate} \otimes \frac{x}{\tate}
  \end{split}
\end{equation}
Clearly, $P(\hM^{(2)})$ is invariant under the
action~\eqref{eq:21} (over $\QQ$). Moreover, $P(\hM^{(2)})$
belongs to the kernel $\curly{I}$ of the multiplication map
\begin{math}
  \CC\otimes_\QQ\CC \to \CC\,.
\end{math}
As a consequence, we have:
\begin{proposition}
  \label{prop:3}
  The ``connection form'' \eqref{eq:22} and the (logarithm of
  the) hermitian fiber metric on the Heisenberg bundle correspond
  to the images of $P(\hM^{(2)})$ under the two projections
  \begin{displaymath}
    \curly{I} \lto \curly{I}/\curly{I}^2 =
    \Omega^1_{\CC/\QQ}
  \end{displaymath}
  and
  \begin{displaymath}
    \curly{I}\subset  \CC\otimes_\QQ \CC
    \lto \RR(1)\,,
  \end{displaymath}
  respectively.
\end{proposition}
\begin{proof}
  The images under the two projections are, respectively, equal
  to
  \begin{displaymath}
    -d \Bigl(\frac{z}{(\tate)^2}\Bigr)
    + \frac{x}{\tate}\, d \Bigl(\frac{y}{\tate}\Bigr)
  \end{displaymath}
  and
  \begin{displaymath}
    \frac{1}{(\tate)^2} \bigl(\pi_1(z)
    -\pi_1(x)\,\pi_0(y)\bigr)\,. 
  \end{displaymath}
\end{proof}

\subsection{The extension class}
\label{sec:extension-class}
The big period can be obtained as a symmetrization of an
extension class of MHS. Indeed, the weight $-2$ subspace
$W_{-2}\hM^{(2)}\iso \hM^{(1)} \otimes \tate \coin \hM^{(1)} (1)$
is itself a MHS (twisted by $\tate$) defined by
\begin{equation}
  \label{eq:32}
  A^{(1)} = 
  \begin{pmatrix}
    1 & \\ y & \tate
  \end{pmatrix}\,.
\end{equation}
(The data are as for $\hM^{(2)}$, replacing $2$ by $1$.) We thus
have an extension of MHS:
\begin{equation}
  \label{eq:33}
  0 \lto \hM^{(1)}(1) \lto \hM^{(2)}
  \lto \ZZ(0) \lto 0\,.
\end{equation}
Following the procedure explained in ref.~\cite{MR1148383}, it is
seen that the class of the extension~\eqref{eq:33} belongs to
\begin{displaymath}
  \hM^{(1)}_\CC(1) /\hM^{(1)}_\QQ(1)\,,
\end{displaymath}
and it is given by the vector
\begin{equation}
  \label{eq:34}
  e = -\frac{x}{\tate}\, v_1
  - \frac{z -xy}{(\tate)^2} \, v_2
\end{equation}
taken modulo $\hM^{(1)}_\QQ$. This computation can be refined by
noticing (\cite{MR1148383}) that $\hM^{(1)}$ is itself an
extension,
\begin{displaymath}
  0 \lto \ZZ(1) \lto \hM^{(1)}
  \lto \ZZ(0) \lto 0
\end{displaymath}
mapping (over \QQ) to the ``universal extension'' $\hH^{(1)}$:
\begin{equation}
  \label{eq:35}
  0 \lto \QQ(1) \lto \CC
  \lto \CC^\unit\otimes\QQ \lto 0
\end{equation}
obtained by tensoring the standard exponential sequence by
\QQ. Over the complex numbers, we have
\begin{displaymath}
  0 \lto \CC \lto \CC\otimes_\QQ \CC
  \lto \CC^\unit\otimes_\ZZ\CC \iso 
  \CC/\QQ(1)\otimes_\QQ\CC\lto 0\,,
\end{displaymath}
Here we have $\hH^{(1)}_\QQ = \CC$ and 
\begin{math}
  \hH^{(1)}_\CC = \CC\otimes_\QQ\CC\,.
\end{math}
According to the same principle the class of the
extension~\eqref{eq:35} lives in 
\begin{equation}
  \label{eq:36}
  \hH^{(1)}_\CC /\hH^{(1)}_\QQ \iso
  \CC\otimes_\QQ\CC / \CC
  \iso \CC\otimes_\ZZ\CC^\unit\,.
\end{equation}
The image of \eqref{eq:34} in 
\begin{math}
  \CC\otimes_\QQ \CC
\end{math}
is given by
\begin{equation}
  \label{eq:37}
  \Tilde e = 
  - y\otimes x - \tate\otimes \frac{z-xy}{\tate}\,.
\end{equation}
Taking \eqref{eq:37} modulo $\hH^{(1)}_\QQ\iso \CC$ we finally
have
\begin{equation}
  \label{eq:38}
  (\Id\otimes\exp) (\Tilde e) =
  y\otimes e^{-x} + \tate\otimes e^{-(z-xy)/\tate}\,.
\end{equation}
This is the (image of) the class of the extension~\eqref{eq:33}
as computed in ref.~\cite{MR1148383}. It is easily seen that the
element \eqref{eq:38} is invariant under the
transformations~\eqref{eq:21}.
\begin{lemma}
  There is a unique well defined lift of the class \eqref{eq:38}
  to
  \begin{math}
    F^0\hH^{(1)}_\CC = \ker (m\colon \CC\otimes_\QQ\CC \to
    \CC)\,. 
  \end{math}
  This can be obtained by adding to~\eqref{eq:37} a (necessarily
  unique, see ref.~\cite{MR1148383}) element from
  $\hH^{(1)}_\QQ\iso \CC$ to~\eqref{eq:37}. The lift is 
  \begin{displaymath}
    \tate\otimes \tate \cdot P(\hM^{(2)})\,.
  \end{displaymath}
\end{lemma}
\begin{proof}
  We can identify $\hH^{(1)}_\QQ\iso \CC$ inside $\hH^{(1)}_\CC$
  via $a \mapsto a\otimes\tate$. Thus add any such element to
  $\Tilde e$ and consider the image under the multiplication map:
  \begin{displaymath}
    m( \Tilde e + a\otimes \tate) = -z +\tate a\,.
  \end{displaymath}
  It is equal to zero iff $a = z/\tate$, hence
  \begin{align*}
    \Tilde{\Tilde e} &= \Tilde e + \frac{z}{\tate}\otimes \tate \\
    &= - y\otimes x + \tate\otimes \frac{xy}{\tate}
    +\frac{z}{\tate}\otimes \tate
    -\tate \otimes\frac{z}{\tate}
  \end{align*}
  is the required element.
\end{proof}

\bibliography{general} \bibliographystyle{hamsplain}

\def\polhk#1{\setbox0=\hbox{#1}{\ooalign{\hidewidth
  \lower1.5ex\hbox{`}\hidewidth\crcr\unhbox0}}}
\providecommand{\bysame}{\leavevmode\hbox to3em{\hrulefill}\thinspace}
\begin{thebibliography}{10}

\bibitem{math.CV/0211055}
Ettore Aldrovandi, \emph{On hermitian-holomorphic classes related to
  uniformization, the dilogarithm and the liouville action},
  \mbox{arXiv:math.CV/0211055}, To appear in Commun. Math. Phys.

\bibitem{MR86h:11103}
A.~A. Be{\u\i}linson, \emph{Higher regulators and values of {$L$}-functions},
  Current problems in mathematics, Vol. 24, Itogi Nauki i Tekhniki, Akad. Nauk
  SSSR Vsesoyuz. Inst. Nauchn. i Tekhn. Inform., Moscow, 1984, pp.~181--238.

\bibitem{bei:hodge_coho}
\bysame, \emph{Notes on absolute {H}odge cohomology}, Applications of algebraic
  $K$-theory to algebraic geometry and number theory, Part I, II (Boulder,
  Colo., 1983), Amer. Math. Soc., Providence, RI, 1986, pp.~35--68.

\bibitem{MR95a:19008}
A.~A. Be{\u\i}linson and P.~Deligne, \emph{Interpr\'etation motivique de la
  conjecture de {Z}agier reliant polylogarithmes et r\'egulateurs}, Motives
  (Seattle, WA, 1991), Proc. Sympos. Pure Math., vol.~55, Amer. Math. Soc.,
  Providence, RI, 1994, pp.~97--121.

\bibitem{bloch:dilog_lie}
Spencer Bloch, \emph{The dilogarithm and extensions of {L}ie algebras},
  Algebraic $K$-theory, Evanston 1980 (Proc. Conf., Northwestern Univ.,
  Evanston, Ill., 1980), Springer, Berlin, 1981, pp.~1--23.

\bibitem{MR1148383}
\bysame, \emph{Function theory of polylogarithms}, Structural properties of
  polylogarithms, Math. Surveys Monogr., vol.~37, Amer. Math. Soc., Providence,
  RI, 1991, pp.~275--285.

\bibitem{MR95m:18006}
Lawrence Breen, \emph{On the classification of {$2$}-gerbes and {$2$}-stacks},
  Ast\'erisque (1994), no.~225, 160.

\bibitem{brymcl:deg4_I}
J.-L. Brylinski and D.~A. McLaughlin, \emph{The geometry of degree-four
  characteristic classes and of line bundles on loop spaces. {I}}, Duke Math.
  J. \textbf{75} (1994), no.~3, 603--638.

\bibitem{brymcl:deg4_II}
\bysame, \emph{The geometry of degree-$4$ characteristic classes and of line
  bundles on loop spaces. {I}{I}}, Duke Math. J. \textbf{83} (1996), no.~1,
  105--139.

\bibitem{bry:loop}
Jean-Luc Brylinski, \emph{Loop spaces, characteristic classes and geometric
  quantization}, Birkh\"auser Boston Inc., Boston, MA, 1993.

\bibitem{bry:quillen}
\bysame, \emph{Geometric construction of {Q}uillen line bundles}, Advances in
  geometry, Birkh\"auser Boston, Boston, MA, 1999, pp.~107--146.

\bibitem{MR97d:32041}
Jean-Luc Brylinski and Dennis McLaughlin, \emph{Holomorphic quantization and
  unitary representations of the {T}eichm\"uller group}, Lie theory and
  geometry, Progr. Math., vol. 123, Birkh\"auser Boston, Boston, MA, 1994,
  pp.~21--64.

\bibitem{del:symbole}
P.~Deligne, \emph{Le symbole mod\'er\'e}, Inst. Hautes \'Etudes Sci. Publ.
  Math. (1991), no.~73, 147--181.

\bibitem{esn:char}
H{\'e}l{\`e}ne Esnault, \emph{Characteristic classes of flat bundles}, Topology
  \textbf{27} (1988), no.~3, 323--352.

\bibitem{esn-vie:del}
H{\'e}l{\`e}ne Esnault and Eckart Viehweg, \emph{Deligne-{B}e\u\i linson
  cohomology}, Be\u\i linson's conjectures on special values of $L$-functions,
  Academic Press, Boston, MA, 1988, pp.~43--91.

\bibitem{MR49:8992}
Jean Giraud, \emph{Cohomologie non ab\'elienne}, Springer-Verlag, Berlin, 1971,
  Die Grundlehren der mathematischen Wissenschaften, Band 179.

\bibitem{MR1978709}
A.~B. Goncharov, \emph{Explicit regulator maps on polylogarithmic motivic
  complexes}, Motives, polylogarithms and Hodge theory, Part I (Irvine, CA,
  1998), Int. Press Lect. Ser., vol.~3, Int. Press, Somerville, MA, 2002,
  pp.~245--276.

\bibitem{MR99i:19004}
Alexander Goncharov, \emph{Volumes of hyperbolic manifolds and mixed {T}ate
  motives}, J. Amer. Math. Soc. \textbf{12} (1999), no.~2, 569--618,
  \mbox{arXiv:alg-geom/9601021}.

\bibitem{gh:alg_geom}
Phillip Griffiths and Joseph Harris, \emph{Principles of algebraic geometry},
  Wiley-Interscience [John Wiley \& Sons], New York, 1978, Pure and Applied
  Mathematics.

\bibitem{MR94k:19002}
Richard~M. Hain, \emph{Classical polylogarithms}, Motives (Seattle, WA, 1991),
  Proc. Sympos. Pure Math., vol.~55, Amer. Math. Soc., Providence, RI, 1994,
  pp.~3--42.

\bibitem{lang:arakelov}
Serge Lang, \emph{Introduction to {A}rakelov theory}, Springer-Verlag, New
  York, 1988.

\bibitem{rama:reg_hei}
Dinakar Ramakrishnan, \emph{A regulator for curves via the {H}eisenberg group},
  Bull. Amer. Math. Soc. (N.S.) \textbf{5} (1981), no.~2, 191--195.

\end{thebibliography}

\end{document}